\documentclass[11pt]{article}
\usepackage{natbib}
\usepackage[asymmetric]{geometry}
\usepackage{breqn}
\usepackage{stmaryrd}
\usepackage{url}
\usepackage{float,here}
\usepackage{amsthm,bm}
\usepackage{graphicx} 
\usepackage{amsmath, amsthm, amssymb, amsfonts, mathrsfs, mathtools}
\usepackage{epstopdf}
\usepackage{makeidx}
\usepackage{xcolor}
\usepackage{subcaption}
\usepackage{ulem}    
\usepackage{algpseudocode} 
\usepackage[]{algorithm2e}
\usepackage{authblk}
\newcommand\norm[1]{\lVert#1\rVert}

\theoremstyle{definition}
\newtheorem{remark}{Remark}[section]  

\title{Jacobian-Free Newton-Krylov with a globalization method for \\
solving groundwater flow models of multi-layer aquifer systems}

\author{
Raghav Singhal\thanks{Corresponding author: Department of Computer Science, University of California, Davis, CA 95616 USA, \\ \textit{rsinghal@ucdavis.edu}}
\quad 
Emin Can Dogrul \thanks{California Department of Water Resources, Modeling Support Office, 1516 9th Street, Sacramento, CA 95814, USA, \textit{can.dogrul@water.ca.gov}} 
\quad
Zhaojun Bai\thanks{
Department of Computer Science and Department of Mathematics, University of California, Davis, CA 95616 USA, \textit{zbai@ucdavis.edu}
}
}
\date{\today}
\begin{document}
\maketitle
\begin{abstract}
    A Jacobian-free Newton–Krylov (JFNK) method with a globalization scheme is introduced to solve large and complex nonlinear systems of equations that arise in groundwater flow models of multi-layer aquifer systems. We explore the advantages of the JFNK method relative to the Newton-Krylov (NK) method and identify the circumstances in which the JFNK method demonstrates computing efficiency. We perform the validation and efficiency of the JFNK method on various test cases involving an unconfined single-layer aquifer and a two-layer aquifer with both confined and unconfined conditions. The results are validated by the NK method. The JFNK method is incorporated in Integrated Water Flow Model (IWFM), an integrated hydrologic model developed and maintained by California Department of Water Resources. We examine the determinacy of the JFNK's adaptability on practical models such as the California Central Valley Groundwater-Surface Water Simulation Model (C2VSim).
\end{abstract}


\section{Introduction}

\paragraph{Background and motivation.} Hydrological processes enclose the continuous cyclical movement of water through the Earth's atmosphere, surface, and subsurface. This process involves the movement of water through a number of stages, including precipitation, infiltration, evaporation, transpiration, condensation, runoff, and groundwater flow. Groundwater flow is an essential component of the hydrological cycle. The water enters the ground, travels through soil and rock, and then through aquifers. Gravitational forces and pressure differences propel it as it moves through the aquifers, gradually flowing through the pore spaces in the aquifer material. Finally, it returns to the surface water bodies. This movement is crucial for sustaining water availability, preserving ecosystem health, and effectively managing water resources. The quantitative assessment of groundwater flow has been a significant concern for engineers, hydrologists, and applied mathematicians for over a century, and the methodologies for its characterization are among the primary concerns in groundwater hydrology.

\paragraph{Mathematical modelling and analytical approaches.} We know that we cannot mathematically model all aspects of real-world movement of groundwater. It is of the utmost significance that a modeling adheres to the fundamental principles of conserving mass, momentum, and energy. One can use partial differential equations (PDEs) to provide a mathematical description of groundwater flow, as they govern the movement of water through porous media (such as soil or rock). These equations take into account the physical processes involved in groundwater flow. These processes include fluid dynamics, porosity, and permeability. The predominant governing equation for groundwater flow is the combination of Darcy's Law \cite{darcy1856fontaines} and the continuity equation, yielding the groundwater flow equation, which is formulated mathematically and analyzed by Boussinesq \cite{boussinesq1904recherches} in 1904. The nonlinear Boussinesq equation (BE) approximates the groundwater flow in an unconfined aquifer. The resulting equation describes the spatial and temporal changes in the hydraulic head inside the aquifer. These equations take into consideration things like specific storage, which means how much water one can store in the aquifer, and hydraulic conductivity, which determines the ease of water movement across the medium.

Boussinesq equation is also extensively used in various domains of hydrology, such as irrigation and drainage \cite{dan2012boussinesq}, baseflow investigations and recession analysis \cite{rupp2006use} and constructed subsurface wetlands \cite{lockington2000similarity}, and catchment hydrology \cite{troch2013importance}. Solutions to BE are significant since they give important insights into the water table's reaction to changes in stream levels and help with the quantification of the exchange between the aquifer and stream. Moutsopoulos and Konstantinos \cite{moutsopoulos2013solutions} applied the domian decomposition method to obtain a semi-analytical solution of the BE with a nonlinear Robin boundary condition. Bartlett and Porporato \cite{bartlett2018class} came up with a set of analytical solutions to BE for both horizontal and sloping aquifers. These solutions take into account sources and sink factors, such as water seeping through the bedrock. The solutions have an explicit dependency on spatial coordinates and are restricted by particular time-dependent boundary conditions and space-dependent initial conditions derived from the general forms of the true solutions. Hayek \cite{hayek2019accurate} introduced a novel approach to addressing BE in horizontal unconfined aquifers resulting from sudden changes in the boundary head by employing an empirical function characterized by four parameters, which may be determined through a numerical fitting procedure. 

\paragraph{Numerical approaches.} The literature reveals analytical solutions to BE for only a limited number of circumstances with highly simplified initial and boundary conditions. Determining the analytical solution is generally challenging due to its non-linear nature, particularly when considering the flow in multiple heterogeneous aquifer layers with complex initial and boundary conditions and under complex sources and sinks such as recharge and pumping. Discretizing the governing PDEs using any traditional method such as finite difference (FD), finite volume (FV), or finite element (FE) method, leads to a nonlinear system of algebraic equations. It poses enormous challenges in terms of computational efficiency and memory requirements. A common approach to solve these nonlinear system of equations is the Newton method, also known as the Newton-Raphson method. The Newton method needs to forming and storing the Jacobian matrix at each time step, which can be expensive. Since the resulting linear system does not have an exact solution, it is generally necessary to adopt an iterative solution method. Thus, the resulting method is known as the inexact Newton method.  Employing a direct method to solve a large system of linear equations may prove inefficient. The Krylov subspace method is the most prominent iterative method for solving large linear systems. Thus, the integration of the Newton method for solving nonlinear problems with a Krylov subspace based iterative approach is referred to as the Newton-Krylov method. In Newton's method one can also approximate the Jacobian by using difference approximation, which is usually referred to as a Jacobian-free approach \cite{kelley1995iterative, maxwell2013terrain}. In contrast, Jacobian-Free Newton-Krylov (JFNK) method eliminates the need for explicit computation of the Jacobian matrix, where we approximate the Jacobian-vector product using the evaluation of nonlinear function values, which is the focus of this study. JFNK method is an exemplary form of a nested iteration approach that involves a minimum of two and usually four levels. The nomenclature of this approach is derived from its fundamental components, namely the iteration over the Newton step and the iteration for constructing the Krylov subspace, from which each Newton correction is obtained. A globalization method is frequently needed when working outside of the Newton iterations. A preconditioner is commonly required to accelerate the convergence in the Krylov iterations for solving the linear system. It can be implemented as either a direct or iterative procedure. The common approaches for effective preconditioners are incomplete lower-upper (ILU) factorization, Newton–Krylov–Schwarz (NKS) \cite{cai1994newton}, and the multigrid method \cite{kothari2023multigrid, mavriplis2002assessment}; however, physics-based preconditioners \cite{kadioglu2017second, kadioglu2013jacobian, kadioglu2023Jacobian} are the most popular. The physics-based preconditioners depend on the underlying physics and numerical discretization. The preconditioning is beyond the scope of this work. However, exploring the use and benefits of Anderson acceleration \cite{walker2011anderson} as an alternative to preconditioning could effectively reduce the need for knowledge of numerical discretization.

\paragraph{Related work on JFNK.} Knoll et al. \cite{knoll2005Jacobian} presented the idea of applying semi-implicit methods as preconditioners to JFNK methods for the accurate time-integration of stiff wave systems. This approach was examined on two-dimensional problems that came under the umbrella of geophysical fluid dynamics and magnetohydrodynamics. Evans and his coworkers \cite{evans2006development} developed a JFNK solution methodology for phase change convection in two dimensions using the incompressible Navier-Stokes equation set and enthalpy as the energy conservation variable. The JFNK method has also been used in a variety of research interests to solve problems such as radiation hydrodynamics, incompressible flow, Fisher equations, Richard's equation, multilevel non-local thermodynamical equilibrium radiative transfer problems, and others (\cite{amir2021jacobian, arramy2024Jacobian, brown1990hybrid, evans2007enhanced, gomez2023Jacobian, jones2001newton, knoll2004Jacobian, kothari2023multigrid, liu2023modified, maxwell2013terrain, parand2017numerical}).

\paragraph{IWFM related work.} To solve the groundwater flow equation in multiple heterogeneous aquifer layers, several numerical models have been developed. PARallel FLOW (ParFlow) \cite{maxwell2025parflow} is a comprehensive hydrological model that simulates both surface and subsurface flows.  ParFlow is a three-dimensional numerical code that uses a Newton-Krylov nonlinear solver and mulitgrid-preconditioned conjugate gradient solver to compute saturated and variably saturated subsurface flow in heterogeneous porous media. HydroGeoSphere \cite{brunner2012hydrogeosphere} simulates the variably saturated flow in the aquifer by solving the three dimensional Richards equation. The Penn State Integrated Hydrologic Model (PIHM) \cite{qu2007semidiscrete} uses a semi-discrete finite volume method to solve the groundwater flow equation coupled with conservation equations representing surface flows, overland flow and infiltration. PFLOTRAN \cite{lichtner2015pflotran} is an open-source, highly parallelised simulation code for subsurface flow and reactive transport. It is engineered for high-performance computing (HPC) and applied in environmental and geoscientific modelling applications, such as groundwater flow, geochemical reactions, contaminant transport, and thermal processes in porous media. Modular Ground-Water Flow Model (MODFLOW) \cite{boyce2020one, mcdonald1988modular} is open-source software that is extensively used worldwide for groundwater flow modeling. It was initially developed by the United States Geological Survey (USGS) in the early 1980s and has gone through several updates since then. 

Integrated Water Flow Model (IWFM) \cite{IWFM_2024} is an integrated hydrologic model that has been developed and maintained by California Department of Water Resources to simulate the movement of water through land surface, root zone and groundwater as well as the flow interactions between land surface and subsurface. An important part of the current edition of IWFM is keeping the non-linear features of surface and subsurface flow processes and how they interact with each other.  At its core, IWFM simulates groundwater heads in a multi-layer aquifer system as well as the flows between the layers. Stream and lake flows are also simulated in IWFM. The impact of stream flows and lakes on the groundwater system is modeled by continuously solving the conservation equations for groundwater, streams, and lakes.

\paragraph{Contributions.} The primary goals of this manuscript are: 
\begin{enumerate}
    \item Develop a JFNK approach with a globalization scheme to solve large and complex nonlinear systems of equations that arise in the groundwater flow; such problems pose challenges because of the multiple aquifer layers and the non-smoothness inherent in the models.
    \item  Explore the benefits, and accuracy of the JFNK method over the Newton method, and identify the specific situations where the JFNK method is more efficient. 
    \item Integrate the JFNK approach into IWFM software and study the impact of the algorithm's adaptability on real-life models such as California Central Valley Groundwater$-$Surface Water Simulation Model (C2VSim) \cite{C2VSimFG_2020}.
\end{enumerate}

To the best of our knowledge, this is the first documented attempt to solve groundwater flow equations discretized by finite element method using the JFNK method. Although ParFlow \cite{maxwell2025parflow} refers to an optional Jacobian-free solution approach in its documentation, there are no details provided regarding any benefits or accuracy associated with using this method.

\paragraph{Organization.} The paper is organized in the following manner: In Section \ref{sec:JFNK_algo}, we provide a brief introduction of the JFNK method. Section \ref{sec1:Model} defines the governing equation applied to the groundwater flow model for multi-layer aquifer system. In Section \ref{sec:Numerical examples}, we present the numerical results. In section \ref{sec:Implementation_IWFM}, we integrate the JFNK methodology in IWFM and finally summarize our findings and potential future work in section \ref{sec:Conclusion}.

\section{Numerical method for nonlinear equations}\label{sec:JFNK_algo}

This section explains the mathematical framework for the numerical solution of large scale nonlinear equations defined as:
\begin{equation} \label{eq:nls}
\textnormal{Find} \;\; h^{\star} \in \mathbb{R}^n \; \textnormal{s.t} \;\; F(h^{\star})= 0, 
\end{equation}
where $F:\mathbb{R}^n \rightarrow\mathbb{R}^n$ is a continuously differentiable function.

Newton method is most frequently used to determine the roots of the nonlinear equation \eqref{eq:nls}. It progressively enhances the sequence $ \lbrace h^{k} \rbrace$ via the initial iterate $h^0$ by applying the formula:
\begin{align}
J(h^{k}) \delta h^{k} &= - F(h^{k}) \label{eq:nls4} \\ 
h^{k+1} &= h^{k}+ \delta h^{k},  \hspace{0.5cm} k=0, 1, 2, ... \label{eq:nls4b}
\end{align}
where $J(h^{k}) \in \mathbb{R}^{n \times n}$ is the Jacobian matrix of the function $F(h)$ with respect to $h$ computed at $h^{k}$. $\delta h^{k}$ is the Newton step (also referred as the update or the correction), and $k$ is the iteration index. This method provides rapid convergence when the initial guess is reasonably close to the actual root of $F(h)$ [\cite{wright2006numerical}, Theorem 11.2].

Generally, numerical solutions of the nonlinear system \eqref{eq:nls} begin with an initial iterate $h^0$ and iteratively progress until the following stopping criterion is fulfilled:
 \begin{equation}
\norm{F(h^{k})} \leq \tau_{r} \norm{F(h^0)} + \tau_{a} \label{NSL_eq2}
\end{equation} 
where $\tau_{r}$ and $\tau_{a}$ are relative and absolute tolerance, respectively. Later on, we will refer to the iterations index denoted by $k$ as the nonlinear iterations or the Newton iterations.

 \subsection{Inexact Newton Method}
In certain cases, as is the case in IWFM, the Newton method may encounter a large non-symmetric Jacobian matrix at each nonlinear iteration where we do not have the exact solution to the system of linear equations \eqref{eq:nls4}. To solve a large non-symmetric linear system of equations, one needs to implement an iterative procedure that ensures that the solution converges locally. This iterative approach is commonly known as an inexact Newton method \cite{dembo1982inexact}: 
\begin{equation}
\norm {F (h^{k}) + J(h^{k}) \delta h^{k}} \leq \eta(k) \norm{F (h^{k})}. \label{eq:inexact_ns}
\end{equation}
The Newton step must satisfy \eqref{eq:inexact_ns}. This indicates that the degree of tolerance with which we solve the linear problem \eqref{eq:nls4} on each nonlinear iteration corresponds to the current nonlinear residual. The parameter $\eta(k)$, which represents the forcing term, is a positive constant less than unity. The forcing term can vary dynamically for each nonlinear iteration according to the reduction of the nonlinear function $F$, which will be further discussed in $\S$\ref{subsection:Comments on JFNK Algorithm}. The forcing term governs the local convergence characteristics of the approach. By selecting $\eta (k)$ effectively, one can get the desirable rate of local convergence, up to the accuracy of the exact Newton method (usually quadratic) \cite{dembo1982inexact}.

\subsection{JFNK Method}
Calculating a Newton update, $\delta h^{k}$, in \eqref{eq:nls4} either using the exact or inexact Newton methods constitutes the majority of computational effort. In complex flow models defined on a large dimension, computing Jacobian matrix $(J(h^{k}))$ and solving Newton step ($\delta h^{k}$) in eq. \eqref{eq:nls4} using direct methods like LU decomposition either prohibitive or proves ineffective. This is because the size of the Jacobian matrix is either too big, or the decomposed Jacobian matrix has too many non-zero entries, which are called fill-ins.

The Newton step can be approximated using iterative methods, such as Generalized Minimal RESidual (GMRES) \cite{saad1986gmres} and Bi-Conjugate Gradient STABilized (BiCGSTAB) \cite{van1992bi} methods. Both methods are based on the Krylov subspace approach \cite{saad2003iterative}. To solve the linear system defined in \eqref{eq:nls4}, we use Krylov subspace $K^{p}(J: r_{0})$
\begin{equation}
    K^{p}(J: r_{0}) = \textnormal{span} \lbrace r_{0}, Jr_{0}, J^{2}r_{0}, ..., J^{p-1}r_{0}\rbrace
\end{equation}
where $r_0$ is the initial linear residual, $r_0 = -F(h^{0}) - J \delta h^0$, and $\delta h^0$ represents an initial guess, which is usually taken as zero. The Newton step ($\delta h^{p}$) in \eqref{eq:nls4} is drawn from the subspace spanned by the Krylov vectors $\lbrace r_{0}, Jr_{0}, J^{2}r_{0}, ..., (J)^{p-1}r_{0} \rbrace$, and it can be written as
\begin{equation}
\delta h^{p} = \delta h^{0} + \sum_{l= 0}^{p-1} \phi_{l} J^{l} r_{0} \label{eq5}
\end{equation}
where the nonlinear iteration index $k$ has been omitted. The Krylov iteration index is denoted by $p$. The scalars denoted by $\phi_{l}$ minimize the residual $\norm{F(h^{p}) + J(h^{p}) \delta h^{p}}$. 

We use Arnoldi orthogonalization method to generate an orthonormal basis $q_{1}, q_{2}, q_{3}, ..., q_{p}$ of the Krylov subspace $K^{p}(J: r_{0})$. Therefore,
\begin{equation}
\delta h^{p} = \delta h^{0} + \sum_{l= 0}^{p-1} \beta_{l} q_{l+1} \label{eq6}
\end{equation}
where the coefficients, $\beta_{l}$ minimize the residual. In order to find the unknown $\beta_{l}$'s, we use GMRES \cite{saad1986gmres} method in this work. The iterations used to get the desired $\delta h^{p}$ is referred to as Krylov inner iterations.

Note that using the true Jacobian in evaluating Newton step using the Krylov method is called "Newton-Krylov (NK)" method. However, in Eq.\eqref{eq5}, the Krylov vectors only requires the Jacobian action in the form of matrix-vector products, which can be approximated by
\begin{equation}
J \cdot \delta h \approx \frac{F(h+ \epsilon \delta h)- F(h)}{\epsilon} \label{eq7}
\end{equation}
where $\epsilon$ is a small perturbation. The approximation has an error of $O(\epsilon)$, and we discuss the selection of the perturbation parameter in $\S$ \ref{remark::perturbation_parameter}. Thus, the approximation used in equation \eqref{eq7} yields a Jacobian-free (JF) method to compute the matrix-vector product $J \cdot \delta h$. The combined action of JF and NK is referred to as the JFNK algorithm. Additional details can be found in \cite{knoll2004Jacobian} and more recent work \cite{an2011finite, kadioglu2023Jacobian, liu2020finite}.

\subsection{Globalization}
 The convergence of the Newton method is strongly dependent on the selection of the initial guess in close proximity to the solution. To enhance the likelihood of convergence, the Newton method must be globalized when suitable initial iterates are not present. To illustrate the globalization method, we define the function $f:\mathbb{R}^n \rightarrow \mathbb{R}$ as a smooth function such that $f(h) = \frac{1}{2} F(h) F(h)^{T}$ for $h \in \mathbb{R}^n$.

 Globalization is designed to evaluate if a step yields adequate improvement towards a solution and, if required, change it to get a step that does provide satisfactory advancement. In order to overcome overshooting the value, we use the line search method, where every loop in a line search method calculates a search step length $\lambda$ and then finds the suitable length to go along that direction. The iteration is given by 
 \begin{equation}
     h^{k+1} := h^{k} + \lambda \delta h^{k}
 \end{equation}

 The efficacy of a line search method relies on the appropriate choice of both the direction and the step length. Line search techniques need $\delta h^{k}$ to be in a descent direction, which means that $(\nabla f)^{T} \delta h^{k} = F(h)^{T} J(h) \delta h^{k} < 0$. This condition ensures that the function $f$ may be minimized along this direction. We apply the backtracking line search method. To begin this method, choose the positive initial step $\lambda$ and define the two constants $\alpha$ and $\rho$, both within the open interval $(0, 1)$. Next, repeatedly verify that the following criterion is satisfied:
\begin{equation}
    f(h^{k} + \lambda \delta h^{k}) \leq f(h^{k}) + \alpha \lambda  (\nabla f)^{T} \delta h^{k} \label{eq::back_tracking_ls}
\end{equation}

If the criterion \eqref{eq::back_tracking_ls} is not satisfied, then update the step size ($\lambda$) by multiplying with $\rho$ (i.e., $\lambda \gets \rho \lambda$), and check the criterion again. Continue the procedure until the criterion is not satisfied. 

The starting value of $\lambda$ is set to 1. The contraction factor, represented by the symbol $\rho$, is assigned the default value $1/2$ and $\alpha = 10^{-4}$. The convergence analysis and parameter details of the method can be found in \cite[$\S 3.2$]{wright2006numerical}.  A comprehensive study on the various globalization methods and application to solve fully-coupled solution of the Navier-Stokes equations can be found in \cite{pawlowski2006globalization}. 

\subsection{JFNK with Globalization} 
In this subsection, we present a pseudo-code of the JFNK approach with a globalization method.
\begin{enumerate}
\item Input the initial iterate $h^0$, and line search parameters $\alpha$ and $\rho$.
\item For $k = 0, 1, . . . $, until {$h^k$} converge, DO
\begin{enumerate}
\item choose $\eta(k) \in [0, 1)$, and use a Jacobian-free solver to solve the following linear equations
\begin{equation}
J(h^{k}) \delta h^{k} = - F(h^{k}) \label{alg:eq1}
\end{equation}
to obtain an inexact Newton direction $\delta h^{k}$ such that
\begin{equation}
\norm{F (h^{k}) + J(h^{k}) \delta h^{k}}_{2} \leq \eta(k) \norm{F (h^{k})}_{2} \label{alg:eq2}
\end{equation}

\item Globalization: Set $\lambda = 1$, and obtain the value of $\lambda  \leftarrow \rho \lambda$ which satisfies
\begin{equation}
f(h^{k} + \lambda \delta h^{k} ) \leq f(h^{k}) + \alpha \lambda \nabla f^{T} \delta h^{k}
\end{equation}
where $f$ is a nonlinear function defined by $f(h) = \frac{1}{2}\norm{F(h)}^{2}$ and  $\nabla f(h)^{T} = F(h)^{T} J(h)$.
\item Update: $ h^{k+1} = h^{k} + \lambda \delta h^{k}$. 
\end{enumerate}
\end{enumerate}

In the following $\S$\ref{subsection:Comments on JFNK Algorithm}, we discuss the critical parameters used in this work. 

\subsection{Remarks}\label{subsection:Comments on JFNK Algorithm}
This section provides a concise discussion of the potential issues and the critical parameter values considered in this work. 
\begin{remark}[Stopping criterion]
The stopping criterion defined in \eqref{NSL_eq2} for the outer Newton loop can be advantageous when dealing with parabolic PDEs because the solution may evolve over time and the relative tolerance can adapt to these changes, whereas absolute tolerance specifies a certain amount of acceptable error \cite{kelley1995iterative}. We are addressing a problem characterized by dimensions $L$ and $T$. In example problems we applied the JFNK method for this study, $L$ is measured in feet and $T$ is defined as one day, week, or month, based upon the flow scenario. In addition, the model may have the magnitude of an initial residual $\norm{F(h^{0})}$ of $O(10^{5})$, which could potentially take longer to converge to the acceptable tolerance for the nonlinear system. 

We adopted a more conceptual stopping criterion, defined as the groundwater head difference between the current and preceding iteration, defined as $\norm{\delta h} \leq \tau_{h}$, where $\delta h$ is computed from Eq. \eqref{eq:nls4}. This criterion is considered prior to the implementation of any backtracking Eq. \eqref{eq::back_tracking_ls}, ensuring that convergence is not improperly reported as a result of reduced step sizes.
\end{remark}
\begin{remark}[Forcing term $\eta(k)$]
The authors in \cite{evans2006development, evans2007enhanced} considered the value of $\eta(k) = \eta = 10^{-2}$. Similarly, the author in \cite{kelley1995iterative} uses the values of $0.2$, $10^{-1}$, and $10^{-2}$, with $10^{-1}$ being his preferred choice. Additionally, \cite{kothari2023multigrid} utilizes a value of the forcing term based on the reduction of the residual norm, $\eta(k) = \textnormal{min}(0.5, \norm{F(h^{k})})$. 
     
If small values of $\eta(k)$ are used during the initial phase of Newton iterations, the problem of oversolving may occur, potentially resulting in an accurate linear solution for an inaccurate Newton correction. This could lead to an improper Newton update and a drop in the Newton convergence. Therefore, it is necessary to adapt a loose tolerance in early Newton iterations and progressively tighten it up as the nonlinear iterations progress:
\begin{equation}
 \eta(k)= \left\{\begin{array}{cc}
& \gamma_{ini}, \;\;\;\;\;\;\;\;\;\;\;\;\;\;\;\;\;\;\;\;\;\;\;\;\;\;\;\;\;\;\;\;\;\;\;\; \textnormal{if} \quad   \norm{F(h^{k})}\geq r,\\
& \textnormal{min} \left(\gamma_{ini}, \frac{\norm{F(h^{k})}}{ \norm{F(h^{k-1})}} \right),\;\;\;\;\;\;\;\; \textnormal{if} \quad   \norm{F(h^{k})} < r.
\end{array}\right.
\end{equation} \label{physica_eta}
where $\gamma_{ini} = 0.99$ and $r = 0.625$ are taken from \cite{seinen2018improving}. We have implemented this dynamic forcing term within the framework of our simulation.

\end{remark} 

\begin{remark}[Perturbation parameter $\epsilon$] \label{remark::perturbation_parameter}

The selection of the perturbation parameter $\epsilon$ in the approximation \eqref{eq7} is sensitive to the scaling of $h$ and $\delta h$. $\epsilon$ is determined using established techniques for computing the Jacobian. In this work, we consider the value of $b$ to be $10^{-6}$. By taking the average of all $\epsilon_j$'s and scaling it, we obtain
\begin{equation}
\epsilon = \frac{1}{n\|\delta h\|^2} \left( \sum_{i=1}^n b|h_i| + b \right)  \label{eq16}
\end{equation}
where $n$ represents the dimension of the system \cite{knoll2004Jacobian}.
\end{remark} 
\begin{remark}[Globalization method]
 Few studies \cite{kadioglu2023Jacobian, lemieux2014second} suggest that the line search (LS) method is not required when dealing with parabolic problems, as the initial condition provides a suitable initial guess for the Newton method to converge. However, we have observed that when dealing with complex models or considering larger time steps in the model, the solution may not converge.
 
 We have limited our LS iterations to three per Newton iteration because if the LS decreases the step size by an unreasonably small amount while the Jacobian remains non-singular, the effectiveness of the Newton direction is impaired.
\end{remark}

\section{Model and discretization} \label{sec1:Model}
This section provides a brief description of the governing PDEs and discretization resulting from the groundwater flow model.

\subsection{Groundwater flow model for a multi-layer aquifer system}
Groundwater flow is inherently a three dimensional flow process. However, solving the three-dimensional groundwater flow equation for an unconfined aquifer using finite difference or finite element methods can be costly because $(i)$ it represents a moving boundary problem which requires updating the simulation grid in the vertical direction \cite{neuman1971analysis}, and $(ii)$ the top boundary condition is a highly non-linear equation in the case that the aquifer is unconfined \cite{bear2012modeling}.

\begin{figure}[t]
\centering
\includegraphics[width=0.65\textwidth]{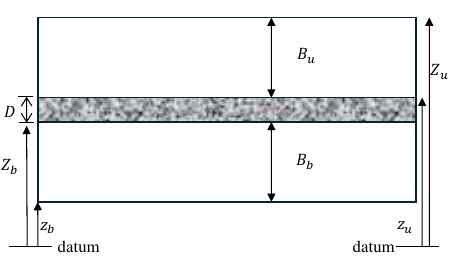} 
\caption{{\sl Present aquifer geometry with $D$ is the thickness of the aquitard. $B_u$ and $B_{b}$ are the aquifer thickness; $Z_u$, $Z_{b}$, and $z_u$, $z_{b}$ are the elevations of the aquifer for the top and bottom layers, respectively.}} \label{fig1}
\end{figure}

To avoid these complications, three dimensional groundwater equation can be integrated across the depth of an aquifer layer using the Dupuit-Forchheimer (DP) approximation which assumes flow within the aquifer layer is primarily horizontal \cite{bear2012modeling}. In an aquifer system where there is substantial vertical flow, the aquifer can be divided into multiple "computational" layers where DP approximation is still applied to each computational layer but the vertical flow can still be simulated. The higher the number of computational layers, the more accurate the simulation of vertical flow becomes. In fact, this approach is implicitly used in many groundwater models including MODFLOW \cite{boyce2020one, mcdonald1988modular} where vertical flow within each computational cell is ignored (i.e. DP approximation is applied at cell level).

Using this approach for two aquifer layers separated by an aquitard (Fig.\ref{fig1}), the flow within each layer can be assumed primarily horizontal and the vertical flow exchange between the layers can be expressed thorough a leakage term. For this aquifer system, the groundwater flow equations in the top and bottom aquifer layers can be expressed as follows: 
\begin{align} \label{iwfm_l2}
\frac{\partial S_{s} h}{\partial t } - \nabla . (\textbf{T} \nabla h) +I \cdot V + I_u \cdot V_u - Q = 0 \quad \forall  x\times y \in \Omega,  t >0 
\end{align}
with  initial conditions
\begin{equation}
    h(x,y, 0) = h_0(x, y), \quad x, y \in \Omega \label{ic}
\end{equation}
and Dirchlet and Neumann boundary conditions
\begin{align}
    h(x,y, t) &= h_b(x, y, t), \quad x, y \in \partial\Omega  \label{bc1}\\
   \textbf{T} \frac{\partial h}{\partial \mathbf{n}} &= q_n(x, y, t), \label{bc2}
\end{align}
where 
$\Omega$ is rectangular domain defined as $[x_{a}, x_{b}] \times [y_{a}, y_{b}]$;
$x$ and $y$ are the space variables in horizontal and vertical directions, respectively ($L$);
$t$ is the time ($T$);
$h$ is the head in the aquifer ($L$); $h_u$ for the top layer ($L$), $h_b$ for the bottom layer ($L$);
$\mathbf{n}$ is the direction normal to the boundary;
$S_{s}$ is the storativity (dimensionless);
$V$ is the vertical flow between an aquifer layer and the layer below $(LT^{-1})$;
$V_u$ is the vertical flow between an aquifer layer and the layer above $(LT^{-1})$;
$I$ and $I_u$ are step functions $(dimensionless)$;
$Q$ is source or sink $(LT^{-1})$;
$D$ is the thickness of the aquitard $(L)$;
$\textbf{T}$ is transmissivity $(L^{2}T^{-1})$:
\begin{equation}
\textbf{T} =  \begin{cases} 
K B & \text{if } h \geq Z \\
K (h - z) & \text{if } h < Z 	\label{iwfm_l6}
\end{cases}
\end{equation}
with $K$ for the horizontal hydraulic conductivity of the aquifer layers ($LT^{-1}$);
$K_v$ for the vertical hydraulic conductivity of the aquitard separating the two aquifer layers $LT^{-1}$;
$z$ for the elevation of the aquifer bottom; $z_u$ for top layer, $z_b$ for bottom layer $(L)$;
$Z$ for the elevation of the aquifer top; $Z_u$ for top layer, $Z_b$ for bottom layer $(L)$;
$B$ and the aquifer thickness; $B_u$ for top layer, $B_b$ for bottom layer $(L)$

The governing equation \eqref{iwfm_l2} is nonlinear in nature and has a nonsmooth solution. When equation \eqref{iwfm_l2} is written for the top aquifer layer, the fourth term drops out since there are no layers above the top aquifer layer; when it is written for the bottom layer, the third term drops out since there are no layers below the bottom layer. To represent this behavior, step functions $I$ and $I_u$ are introduced into equation \eqref{iwfm_l2}. $I_u=0$ and $I=1$ when equation \eqref{iwfm_l2} is written for top layer; $I_u=1$ and $I=0$ when it is written for bottom layer. 

Vertical flows between aquifer layers are defined as:
\begin{equation}
V = -V_u=\frac{ K_v}{D}   \left[ \textnormal{max} (h_u, z_u) -  \textnormal{max} (h_b, Z_b)\right] \label{iwfm_l4}
\end{equation}
Storativity, $S_{s}$, is defined as
\begin{equation}\label{storativity}
 S_{s} = \begin{cases} 
S_o B & \text{if } h \geq Z \\
S_y & \text{if } h < Z 
\end{cases}
\end{equation}
where $S_o$ is the storage coefficient $(L^{-1}$), and $S_y$ is the specific yield (dimensionless). Note that equation \eqref{iwfm_l2}  is a general form of the groundwater equation applicable to both confined and unconfined layers; $S_{s}$ term is assigned the appropriate value based on equation \eqref{storativity} depending on if the layer is confined or unconfined. 

Additionally, in equation \eqref{iwfm_l2} , $Q$ needs to be limited with storage when it is a negative value (i.e. pumping):
\begin{equation}
Q = \begin{cases} 
\min\left(|Q_o|, \frac{M}{\Delta t}\right) & \text{if } Q_o < 0 \\
Q_o & \text{otherwise}	\label{iwfm_l9}
\end{cases}
\end{equation}
where $Q_o$ is the original sink/source term ($LT^{-1}$), $\Delta t$ is the time step ($T$), and $M$ is the water storage available in the aquifer ($L$) described as
\begin{equation}
M = \begin{cases} 
S_y (h - z) & \text{if } h \leq Z \\
S_y B +  S_o B (h - Z) & \text{if } h > Z 	\label{iwfm_l10}
\end{cases}
\end{equation}
For reasons that will be apparent later, we will represent the min function appearing in equation \eqref{iwfm_l9} as follows:
\begin{equation}
\min \left(|Q_o|, \frac{M}{\Delta t} \right) = |Q_o|  - \max\left( |Q_o| - \frac{M}{\Delta t} , 0 \right) \label{iwfm_l11}
\end{equation}
We will also use the following equality later:
\begin{equation}
\max(A, B) = B + \max(A - B, 0)		\label{iwfm_l12}
\end{equation}
To introduce continuity, we use Chen-Harker-Kanzow-Smale Jacobian
smoothing method for the max and the step functions \cite{chen1996class}. Specially the max function can be approximated and smoothed as
\begin{equation}
    \textnormal{max}(h,0) \approx \frac{ h + \sqrt{h^2 + \varepsilon}}{2}
\end{equation}
and the step function ($I$) is approximated and smoothed as \cite{chen1996class}
\begin{equation}
    I(h) \approx \frac{1}{ 1 + e^{- \beta h} } 
\end{equation}
where $\varepsilon$ and $\beta$ is a number indicating the degree of smoothness. In the IWFM setting, the values $\varepsilon = 10^{-4}$ and $\beta = 10$ are considered. In the next subsection $\S$ \ref{subsection:discretization}, we define a numerical method to discretize the governing equation  \eqref{iwfm_l2}.

\subsection{Discretization} \label{subsection:discretization}
The scientific community has encountered persistent challenges in numerically solving PDEs, which are discretized using either a FD, FV, or FE method, depending on whether the focus is on computational or mathematical constraints. Adopting the FE method could be beneficial as it provides a solution for the weak form. Also, groundwater modeling studies frequently employ FE approaches, specifically the Galerkin \cite{gockenbach2006understanding, huyakorn2012computational} and Petrov \cite{brooks1982streamline, fries2004review} FE methods, due to their ability to accurately describe complex boundaries. 

The equation \eqref{iwfm_l2} is first semi-discretized in the spatial domain using the Galerkin FE method. We define the finite element space $V_h\subset H_0^1(\Omega)$ as a space of continuous piecewise linear functions and satisfying the boundary conditions (\ref{bc1})-(\ref{bc2}):
\begin{align}
V_h= \{v:v \in L^2(\Omega), v_x \in L^2(\Omega)  \quad \text{and}\quad v|_{\tau}\in \mathcal{P}_r(\tau) \},
\end{align}
where $\mathcal{P}_r(\tau)$ denotes a piecewise continuous polynomial function of degree $r\in \mathbb{N}$ defined on each quadrilateral or triangular element $\tau$ of the spatial domain $\Omega$, defined by
\begin{align*}
\mathcal{P}_r(\tau)=\text{span}\{x^iy^j: \ 0\leq i,j\leq r\}.
\end{align*}  

The variational formulation is obtained by multiplying the equation (\ref{iwfm_l2}) by the test or shape functions $v\in V_h$ for each aquifer layer and then integrating over the entire domain $\Omega$, which is given by:
\begin{align}\label{weakform}
\nonumber
 &  \text{Find}\ h_h \in V_h \text{such that} \\ 
&\iint\limits_{\Omega}  \left(S_s \frac{\partial h_{h}}{\partial t} + I \cdot V + I_u \cdot V_u \right)v \ d \Omega +  \iint\limits_{\Omega} \mathbf{T} \nabla h_h \cdot \nabla v \ d\Omega- \iint\limits_{\partial \Omega} v q_n \partial \Omega = \left(Q, v\right) \quad  \ \forall\;\; v \in V_h.
\end{align}
The last term of the left side of the equation (\ref{weakform}) is obtained by using the boundary condition (\ref{bc2}) such that
\begin{align*}
\mathbf{T} (\nabla h_h) \cdot  \mathbf{n} =\mathbf{T} \frac{\partial h_h}{\partial  \mathbf{n}}=q_{n}
\end{align*}
The finite element solution $h_h$ to the problem (\ref{weakform}) is approximated by the head values at the discrete node points as
\begin{align}
\label{aprox}
\hat{h}_h(x,y,t)=\sum_{j=N(l-1)+1}^{N_l} h_{h,j}(t) \phi_j(x, y),
\end{align}
where $N$ represents the total number of nodes in an aquifer layer, and $N_l$ is the total number of nodes in all $l$ layers in the aquifer system, respectively. $\phi_j$ are the shape functions that depend on the geometric properties of the finite elements, $h_{h,j}$ are nodal hydraulic head values. 

Substituting (\ref{aprox}) into (\ref{weakform}), we obtain the following equation on each layer:
\begin{align}
\nonumber
\iint\limits_{\Omega} \sum_{j=N(l-1)+1}^{N_l}  S_{s_j} \frac{\partial h_{h,j}}{\partial t} \phi_i \phi_j&+ \left(I \cdot V + I_u \cdot V_u \right)\phi_i \ d \Omega +  \iint \limits_{\Omega}\sum_{j=N(l-1)+1}^{N_l}  \mathbf{T} \nabla (h_{h,j}\phi_j)\cdot \nabla \phi_i \ d \Omega \\  
&- \iint\limits_{\partial \Omega} \sum_{j=N(l-1)+1}^{N_l}  q_n \phi_i  \partial \Omega= (Q,\phi_i), \quad i=1,2 \dots N. \label{e1.8}
\end{align}

By applying the mass lumping method \cite{allen1988numerical, duczek2019mass} to the first term of (\ref{e1.8}), we have
\begin{align}
\nonumber
\iint\limits_{\Omega}  \left( S_{s_i} \frac{\partial h_{h,i}}{\partial t} \phi_i + \left(I \cdot V + I_u \cdot V_u \right)\phi_i \right) d \Omega  +& \sum_{j=N(l-1)+1}^{N_l} \iint \limits_{\Omega}  \mathbf{T}  h_{h,j}  \nabla \phi_{i} \cdot  \nabla \phi_{j} \ d \Omega \\ 
&- \iint\limits_{\partial \Omega}   q_n \phi_i  \partial \Omega= (Q,\phi_i), \quad i=1,2 \dots N. \label{e1.9}
\end{align}

After the equation \eqref{e1.9} is applied to all layers in a multi-layer aquifer system, we obtain a $N \cdot N_l$ system of ordinary differential equations for an aquifer system that is made up of $N_l$ layers with unknown groundwater head values at the $N \cdot N_l$ nodal points. This system is fully discretized on the temporal domain using the Crank-Nicolson method, resulting in the following discrete form:
\begin{align}
& 2 \iint\limits_{\Omega} \frac{(S_{s_i}^{m+1} (h_{h,i}^{m+1}-Z)+S_{s_i}^{m} (h_{h,i}^{m}-Z))}{\Delta t}\phi_i \ d \Omega + \iint\limits_{\Omega} \left(I^{m+1} \cdot V^{m+1} + I_u^{m+1} \cdot V_u^{m+1} \right) \phi_i \ d \Omega \nonumber \\ 
+& \iint\limits_{\Omega} \left(I^{m} \cdot V^{m} + I_u^{m} \cdot V_u^{m} \right)\phi_i \ d \Omega  + \sum_{j=N(l-1)+1}^{N_l} \iint \limits_{\Omega} (  \mathbf{T}^{m+1}  h_{h,j}^{m+1}+  \mathbf{T}^{m}  h_{h,j}^{m} ) \nabla \phi_{i} \cdot  \nabla \phi_{j}d \Omega \nonumber\\
-& \iint\limits_{\partial \Omega}   (q_n^{m+1}+q_n^{m}) \phi_i  \partial \Omega =  (Q^{m+1}+Q^m,\phi_i) \quad \hspace{4.5cm} i=1,2 \dots N \label{eqn:CN_weak_form}
\end{align}
where $m+1$ and $m$ represent the current and previous time levels and $\Delta t$ is the time step. The first integral of equation \eqref{eqn:CN_weak_form} in time discretisation reflects the effort required to change the aquifer conditions. For example, consider an instance in which the aquifer changes from confined to unconfined within a simulation interval of $\Delta t$. More details, including details on the element shape and the specific basis function used in the present simulation, can be found in Chapter $3$ of the IWFM theoretical documentation \cite{IWFM_2024}.

We use JFNK method to solve this system of equations defined in \eqref{eqn:CN_weak_form} for the unknown variable $[h_{h}]_{ N \cdot N_{l} \times 1}$, described in $\S$\ref{sec:JFNK_algo}. This avoids constructing and storing the Jacobian matrix $[J]$ of the size ${ (N \cdot N_{l}) \times (N \cdot N_{l}) }$ on each time step. It is also important to note that JFNK enables the adoption of various space and time discretization techniques, such as higher-order ones. Adopting higher-order methods, where the Jacobian matrix is significantly less sparse, may lead to errors in computing the partial derivative.

\section{Numerical results and discussion}\label{sec:Numerical examples}
This section examines the efficiency and validation of the JFNK method using two groundwater models. The first model is for an unconfined single layer aquifer, and the second model is for a two-layer aquifer with confined and unconfined groundwater flow. We validate our results through the NK method with the exact Jacobian. This section also examines the potential challenges that one may encounter when adopting the JFNK algorithm. All computations are conducted on a PC with a 13th Gen Intel(R) Core(TM) i7-1370P processor and 32 GB of RAM.

\subsection{Test case 1: Groundwater flow model for an unconfined aquifer} \label{subsection:test_case1:GW}
To analyze the JFNK algorithm's correctness and efficiency on solving groundwater equations, we first consider a simple test case with a single unconfined aquifer layer. The model domain is a rectangle defined by $x_{a} = y_{a} = 0$ and $x_{b} = y_{b} = 120000$ ft. The parameters in the model are $S_{y} = 0.25$, $K = 100$ ft$/$day, $z=0$ ft, $Z=500$ ft, and $Q=0$ ft/day. Since this is a single-layer aquifer system, parameters related to vertical flow calculations are not defined. Additionally, because this is an unconfined aquifer, storage coefficient, $S_o$, is also not defined. On the left and right sides of the boundary, $h_b = 50$ ft and $h_b = 400$ ft are taken in \eqref{bc1}, respectively, while the no-flow boundary conditions ($q_n = 0$) are imposed on the top and bottom sides of the domain in \eqref{bc2}. The initial flow is set to $h_0 = 400$ ft in \eqref{ic}.

The simplified model of \eqref{iwfm_l2} has a continuous solution. We use the finite difference method to discretize the simplified model. The grid is generated by vertical and horizontal lines intersecting across the points $(x_i,y_j)$ given by
\begin{equation}
x_i=(i-1)\Delta x, \;\;\;\ y_j=(j-1) \Delta y, \quad \mbox{for $i=1,2,\dots, Nx$ and
$j= 1,2,\dots ,Ny$}. \nonumber
\end{equation}
where $Nx$ and $Ny$ denote the number of node points in the $x$ and $y$ directions. The step length along $x$- and $y$-directions are defined as $\Delta x = \displaystyle{\frac{1}{Nx-1}}$ and $\Delta y = \displaystyle{\frac{1}{Ny-1}}$ respectively. 

We use a resolution of $Nx = Ny = 81$ and the cell widths $\Delta x$ and $\Delta y$ to be $1500$ ft in order to conduct a fair analysis of the computational effectiveness of both approaches. We compute the groundwater head for the time period of $4$ years with a time step $\Delta t = 1$ day. The key variable, such as the tolerance of the outer or nonlinear iterations for NK and JFNK methods are defined as $\tau_{h} =10^{-4}$. The restart parameter $p$, tolerance and the maximum number of restarts $M_k$ for the linear solver GMRES is set to $20$, $10^{-6}$ and $500$, respectively.  

\subsubsection{Validation}

\begin{figure}[!ht]
	\centering
	\begin{subfigure}{0.50\textwidth}		\includegraphics[width=1\linewidth]{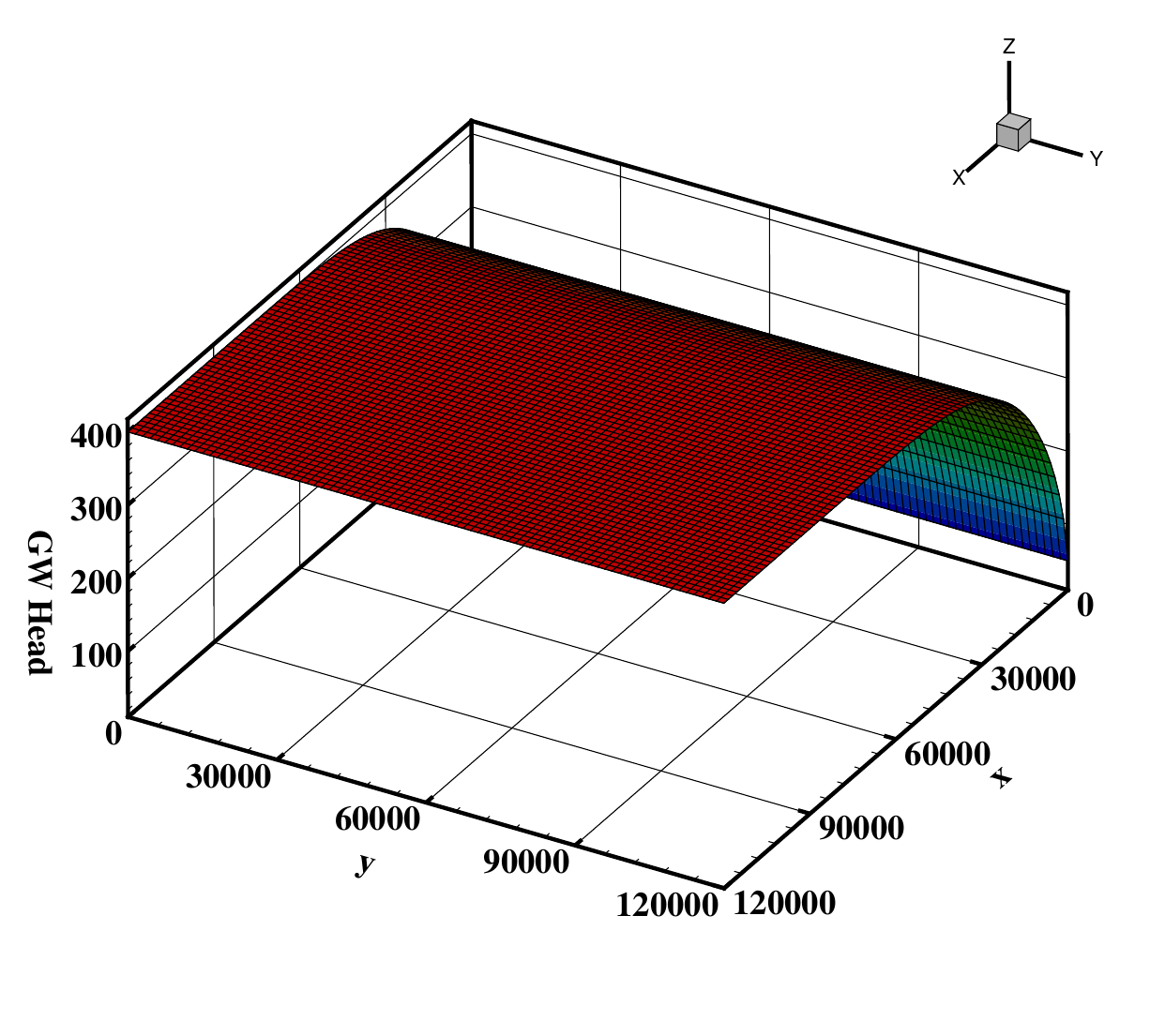}  
		\caption{}
	\end{subfigure}\hfil
	\begin{subfigure}{0.50\textwidth}
	\includegraphics[width=1\linewidth]{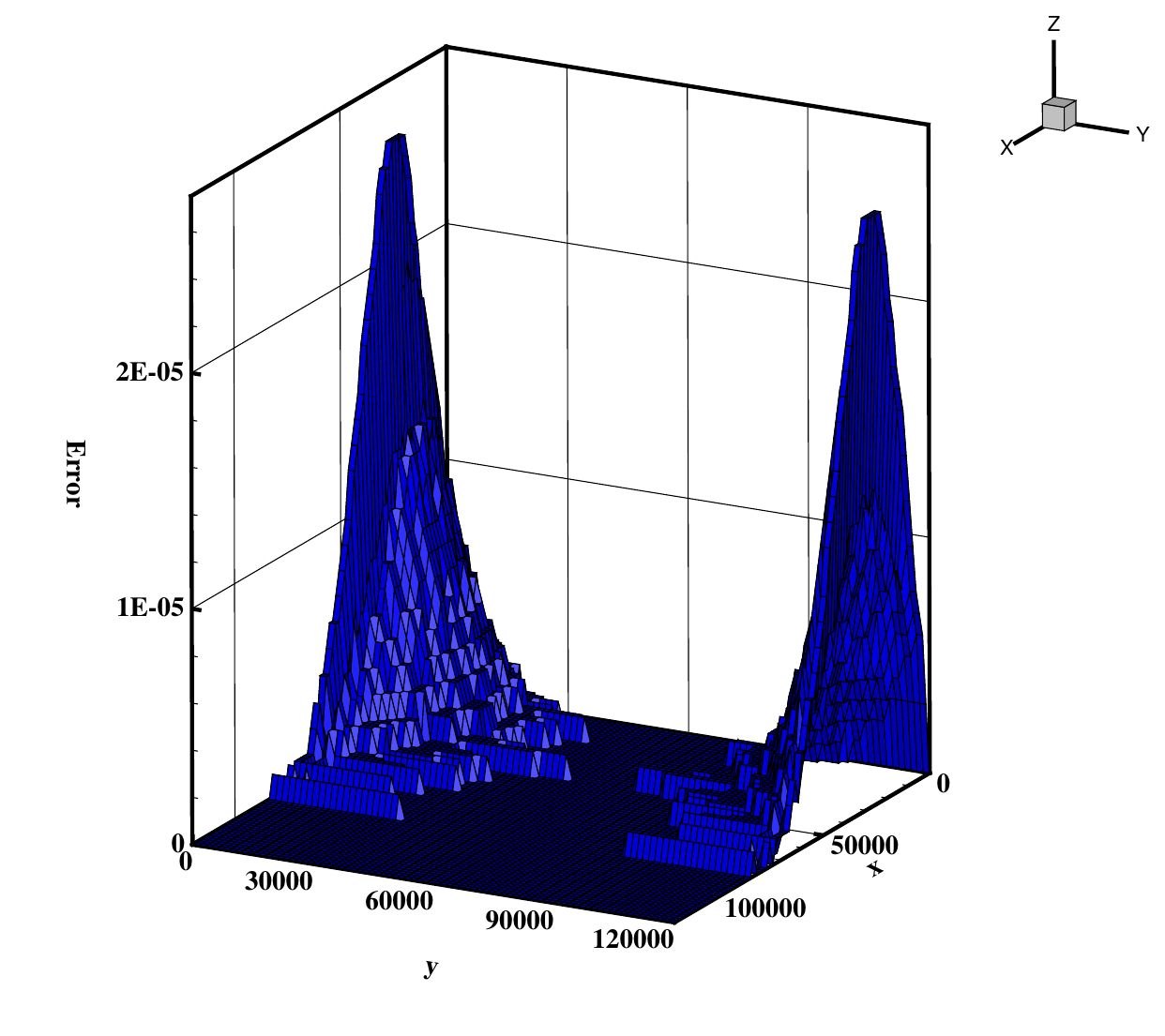}
		\caption{}
	\end{subfigure}\hfil
 \caption{{\sl  Surface plot of (a) the groundwater head using JFNK method and (b) shows the max difference (error) between JFNK and NK solutions at the end of $4$ years for a test case $1$.}} \label{fig1:test_case_1}
\end{figure}

The surface plot of the groundwater head solution from the JFNK method at the end of the four-year simulation period is shown in Figure \ref{fig1:test_case_1}(a). The smoothness of the plot indicates a consistent change in groundwater head throughout the region. 

 Figure \ref{fig1:test_case_1}(b) presents the surface plot exhibiting the error between the results of NK and JFNK. The largest absolute error is at the order of $10^{-4}$, suggesting an adequate correspondence between the results obtained from NK and JFNK.  The surface plot exhibits a spike in the middle of the left and right boundary regions due to the cell width of $1500$ feet, and prevents the finite difference approach from accurately capturing the sudden variation in groundwater between the prescribed boundary conditions in that region. 

\begin{figure}[!ht]
	\centering
	\begin{subfigure}{0.5\textwidth}
		\includegraphics[width=1\linewidth]{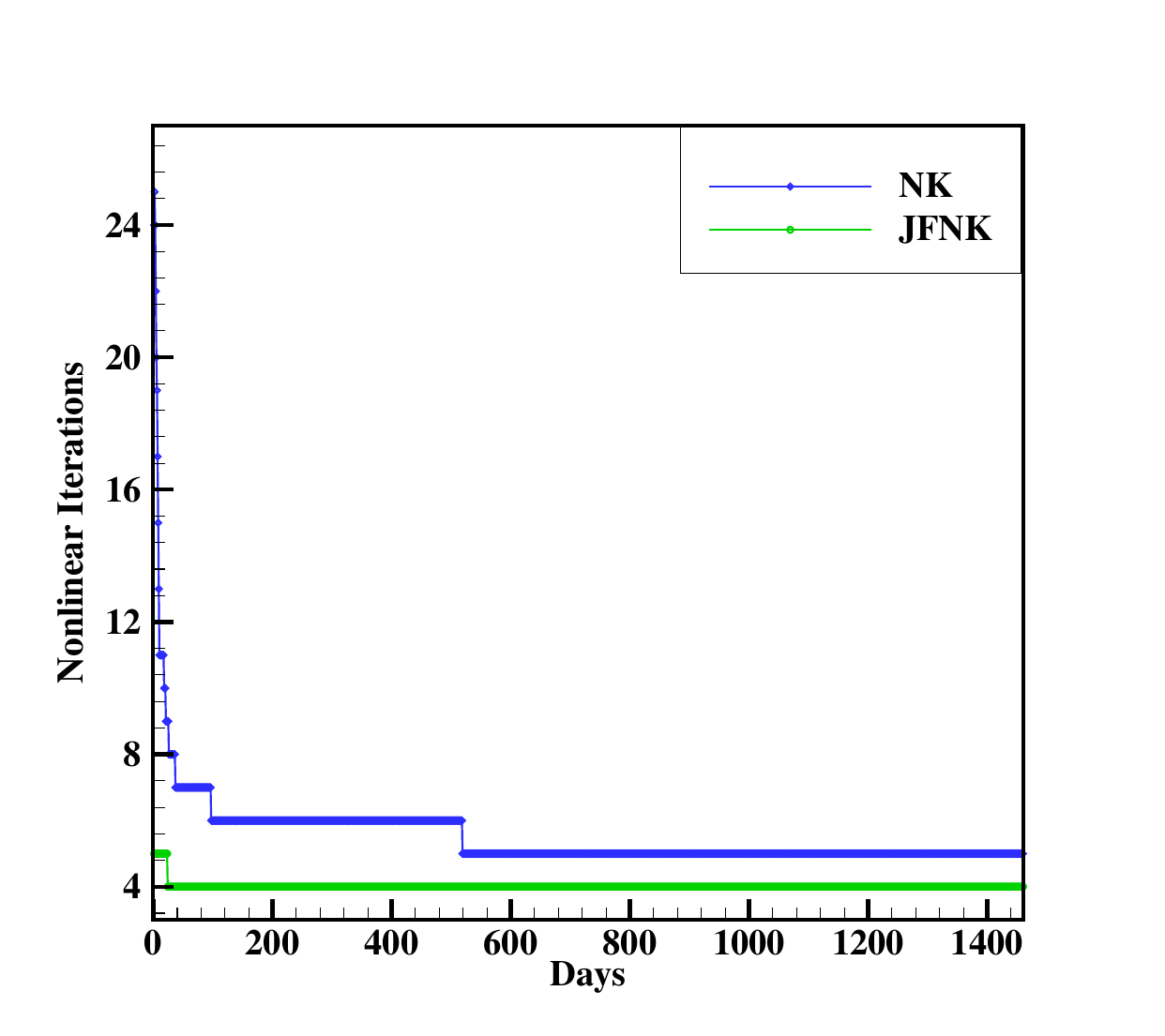}  
		\caption{}
	\end{subfigure}\hfil
	\begin{subfigure}{0.5\textwidth}
		\includegraphics[width=1\linewidth]{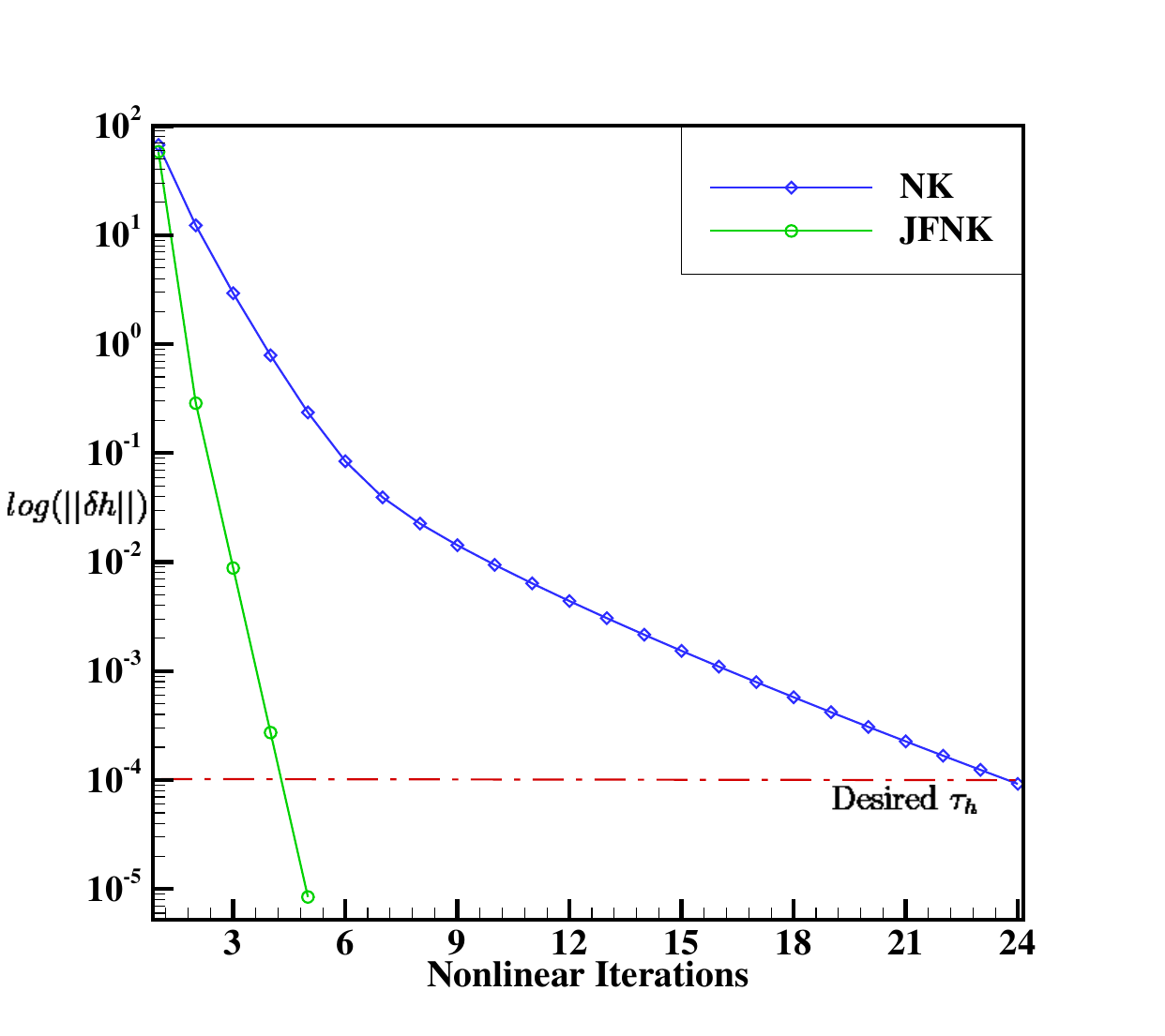}
		\caption{}
	\end{subfigure}
\caption{{\sl (a) The number of nonlinear iterations for each day for both approaches is shown and (b), the decrease in $\norm{\delta h}$ is shown for each nonlinear iteration only in a $T = 1$ one-day simulation for both approaches for test case $1$. } } \label{fig2:test_case_1}
\end{figure}

Figure \ref{fig2:test_case_1}(a) depicts the number of nonlinear iterations corresponding to the NK and JFNK methods over a simulation period of $4$ years, equivalent to $1461$ days. During the initial phase, NK method requires a greater number of nonlinear iterations, but stabilizes around $6$ iterations towards the end of the simulation. However, the JFNK approach requires only $4$-$5$ iterations on each time step throughout the entire simulation period, which is substantially less compared to the NK method. Figure \ref{fig2:test_case_1}(b) shows the decreasing trend of the norm of Newton update ($\delta h$) as the nonlinear iteration progresses, achieving a target tolerance of $\tau_h = 10^{-4}$ for both methods except for the first day where a substantial groundwater head difference between the nonlinear iterations of both approaches is observed in Figure \ref{fig2:test_case_1}(a).

These results demonstrate that the JFNK method accurately captures the head level, providing an alternative method to address the hydrological problem without the necessity of constructing and storing a true Jacobian.

\subsubsection{Efficiency} \label{subsection: eff}
Another performance metric is the CPU time associated with the JFNK method as opposed to the NK method. Also, note that we have not discussed or utilized any preconditioners to solve the linear system \eqref{eq:nls4}. In this subsection, we compare CPU times and discuss elements that impact CPU times.

It is generally computationally expensive to form the Jacobian matrix. The fundamental distinction between the NK and JFNK methodologies is that the NK method creates and stores the Jacobian matrix for each nonlinear iteration. This is expensive and time-consuming, and it can become the bottleneck for the large-scale problems in practice. On the other hand, in JFNK method, the Jacobian matrix is not explicitly formed. Instead, we approximate the Jacobian-vector product by evaluating two functions values (refer to Eq. \eqref{eq7}). Therefore, JFNK method primarily relies on the cost of evaluation of the values of the right-hand side function ($F$ term in equation \eqref{eq:nls4}), namely the RHS function. It is of the utmost importance to be aware of the number of times that the RHS function is called, or the amount of time that it takes to approximate the Jacobian-vector product, to properly gauge the effectiveness of the JFNK method.

Given that the simplified model of \eqref{iwfm_l2} has a smooth solution, employing a preconditioner allows an increase in the convergence rate and reduces the computational effort associated with Krylov vectors. This enhances the practical significance of this experiment. We consider the ILU preconditioner with the GMRES method to solve the linear system \eqref{eq:nls4}. In the NK method, we calculate the exact Jacobian, which allows us to eliminate the implementation of the backtracking line search method, resulting in savings in computational time and reducing the frequency of calls to the RHS function. Note that the GMRES algorithm can employ a maximum of $p \times M_k = 10000$ RHS function calls in a single nonlinear iteration for this test case.

The model runtimes for the above-discussed NK method with a preconditioner and the JFNK method (without preconditioner) are $129.193$ seconds and $107.100$ seconds in MATLAB, respectively. It is evident from the CPU time that JFNK method is more efficient compared to NK method. The profiling data indicates that the total number of calls to RHS function for the whole simulation period of the NK and JFNK methods are $8105$ and $53643$, respectively. This data reveals much more information, including why and when JFNK is superior. It gives a RHS frequency factor ($\varpi$):

\begin{equation}
   \varpi = \frac{\textnormal{Number of RHS calls in JFNK method}}{\textnormal{Number of RHS calls in NK method}} = \frac{53643}{8105} \approx 6.61
\end{equation}
  
  This suggests that if we are using JFNK method, then the frequency of calling the RHS function is approximately \textbf{6.61} times more. Nevertheless, even when we use the NK method with a preconditioner, we obtained better performance with the JFNK method. 
 
 To achieve higher accuracy for the Jacobian-vector product, one may approximate the product in \eqref{eq7} using a second-order finite difference formula. However, this necessitates evaluating $F(h)$ twice, hence increasing $\varpi$ value. By our observations, we found that the adoption of a higher-order approximation of the product did not result in any noticeable improvements in accuracy.

\subsection{Test case 2:  Groundwater flow model for a 2-layer aquifer system with pumping}
In the second test case, we explore groundwater flow in two layers separated by an aquitard of thickness $D$, and have connectivity between the layers (shown in Figure \ref{fig1}). With this model, we aim to test our JFNK results with the IWFM software, which uses the FE method to discretize the model problem \eqref{iwfm_l2} as described in $\S$\ref{subsection:discretization}. 

\begin{figure}[!ht]
\centering
\includegraphics[width=0.6\textwidth]{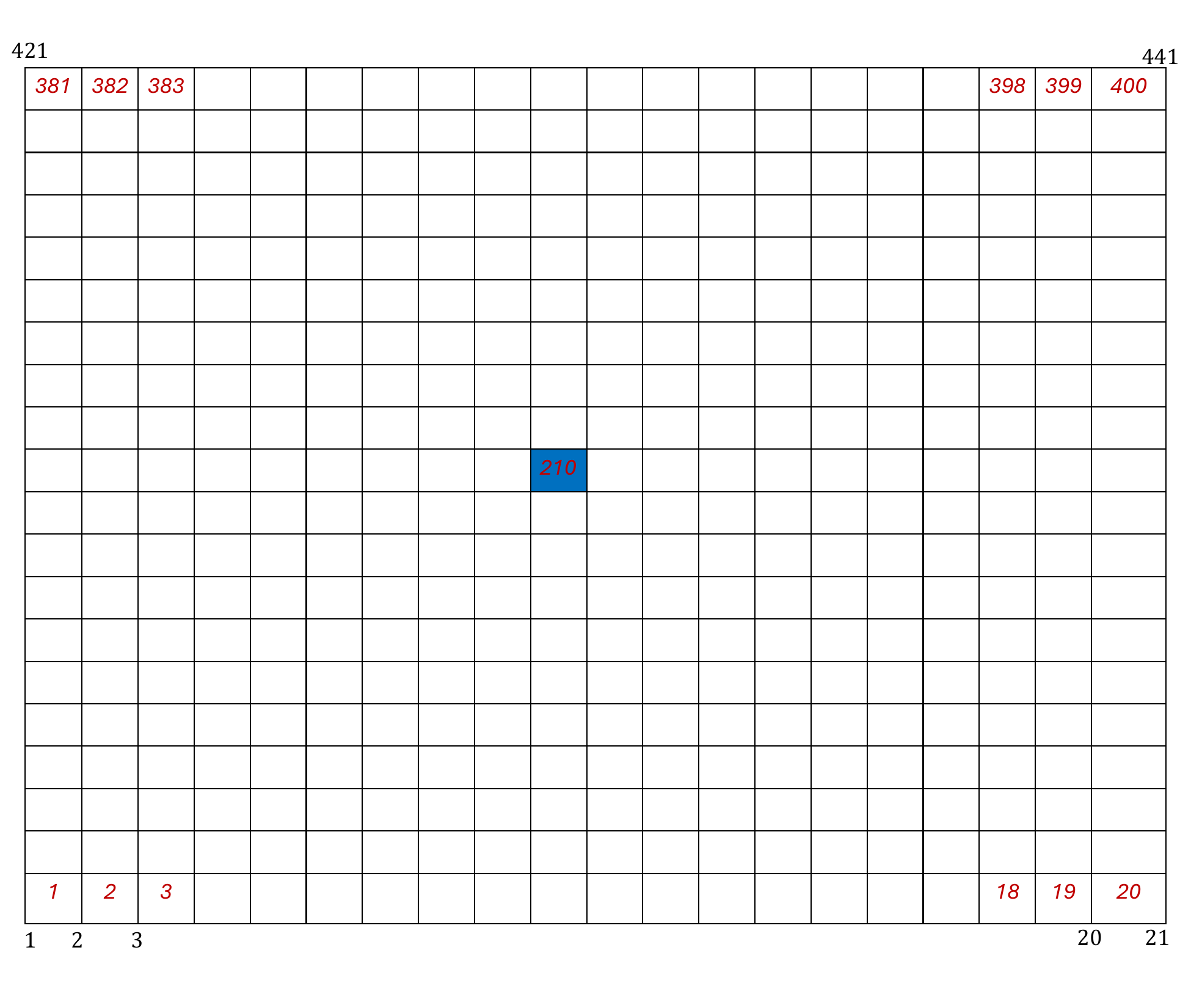} 
\caption{{\sl Schematic representation of a computational domain with a grid size of $21 \times 21$, where nodes are indicated in black and cells in red color. Pumping is defined at $210^{\textnormal{th}}$ cell, which is represented by solid blue color.} }\label{cell_nodes}
\end{figure}

The dimensions of the model domain is the same as the previous test case; i.e. $120000$ ft  by $120000$ ft. We define the number of cells in both directions to be $20$, with a cell width of $\Delta x$ = 6000 ft and $\Delta y$ = 6000 ft for each layer, leading to $882$ node points. The simulation time step is $\Delta t = 1$ day. The other parameters considered are the horizontal hydraulic conductivity $K$ = 100 ft$/$day, specific yield $S_{y}$ = 0.25, and specific storage $S_{0}$$ = 10^{-6}$ 1$/$ft for each layer. The aquitard vertical hydraulic conductivity $K_v=10^{-3}$ ft/day. The ground surface elevation ; i.e. elevation of the top of the upper aquifer layer is $500$ ft, and the aquifer thickness is $300$ ft. The top elevation and the thickness of the bottom aquifer layer are both $170$ ft. The thickness of the aquitard between the two aquifer layers is $D= 30$ ft. A pumping rate, $Q$, of $13068000$  ft$^{3}$/day  is applied at the $210^{\textnormal{th}}$ cell at each aquifer layer (see Figure \ref{cell_nodes}). Although pumping is defined at a cell for this example, IWFM distributes it equally to the surrounding nodes of the cell (nodes $220, 221, 241$ and $242$).  No flow boundary conditions are imposed in the horizontal plane to both aquifer layers. The initial condition $h_0 = 250$ ft is applied to each layer. The tolerance of the nonlinear iterations and other GMRES parameters is the same as in test case $1$ for NK and JFNK methods.

\begin{figure}[!ht]
	\centering
    \begin{subfigure}{0.5\textwidth}
	\includegraphics[width=1\linewidth]{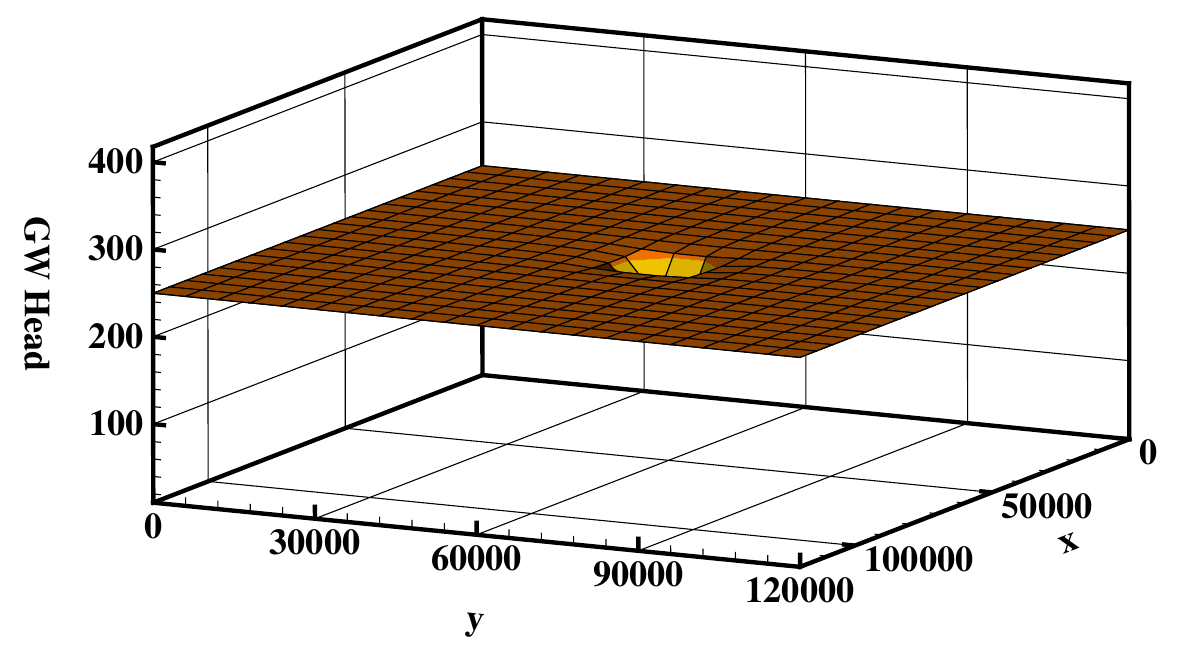}  
	\end{subfigure}\hfil
	\begin{subfigure}{0.5\textwidth}
	\includegraphics[width=1\linewidth]{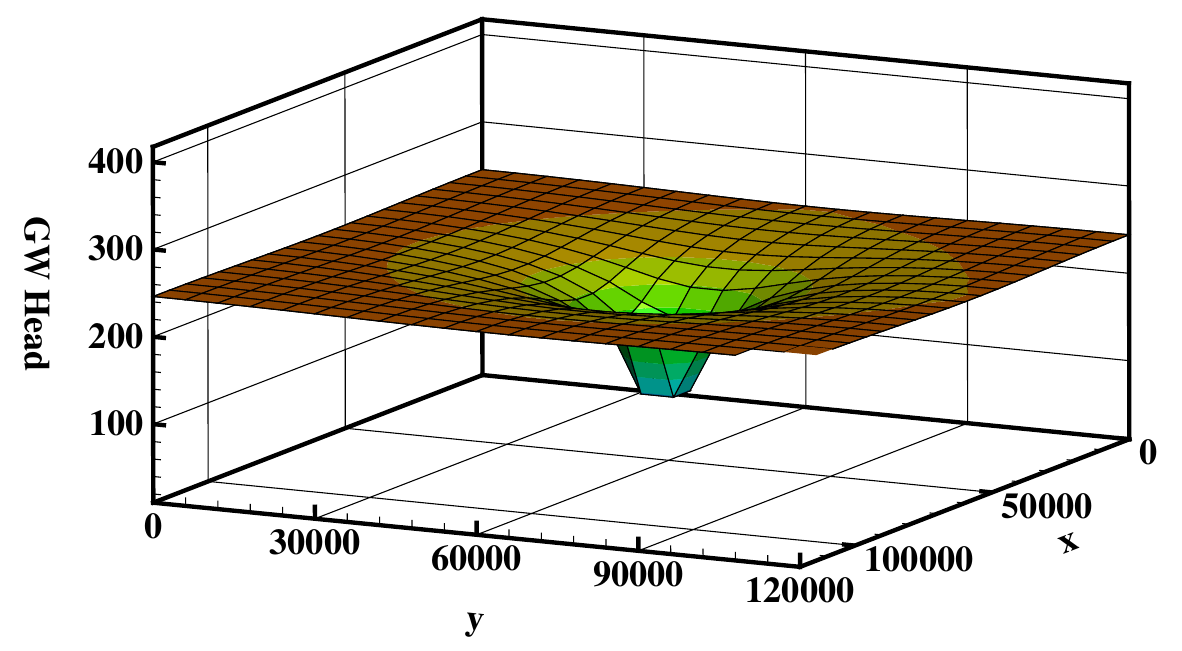}
    \end{subfigure}\hfil
    \begin{subfigure}{0.5\textwidth}
	\includegraphics[width=1\linewidth]{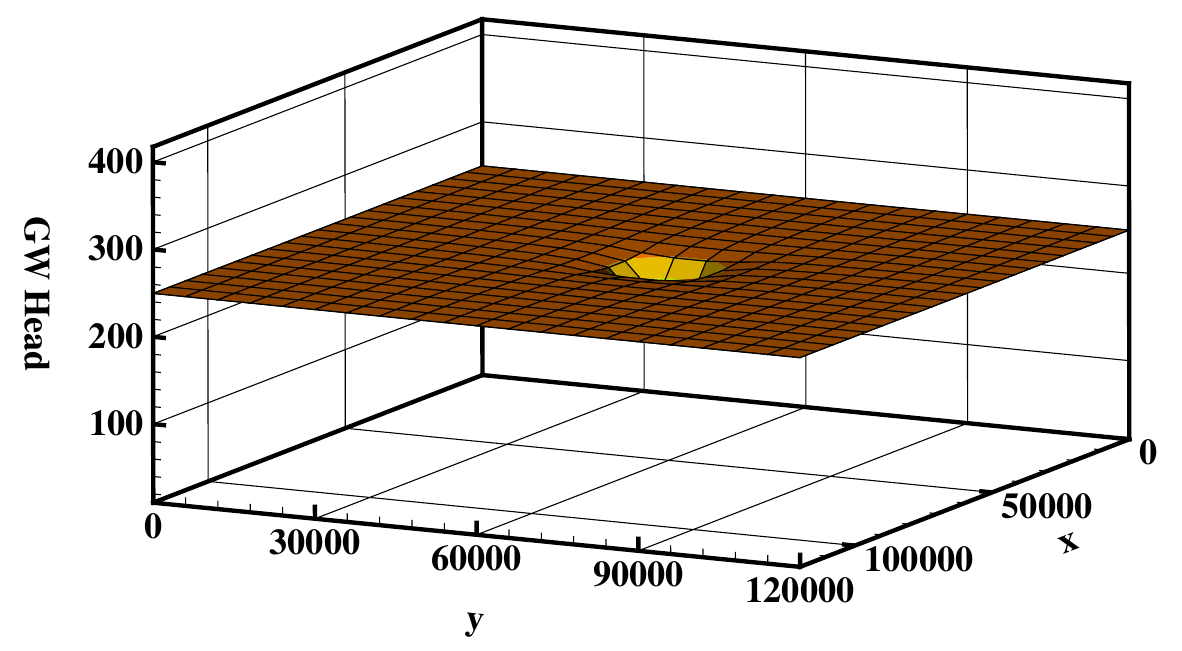}  
	\end{subfigure}\hfil
	\begin{subfigure}{0.5\textwidth}
	\includegraphics[width=1\linewidth]{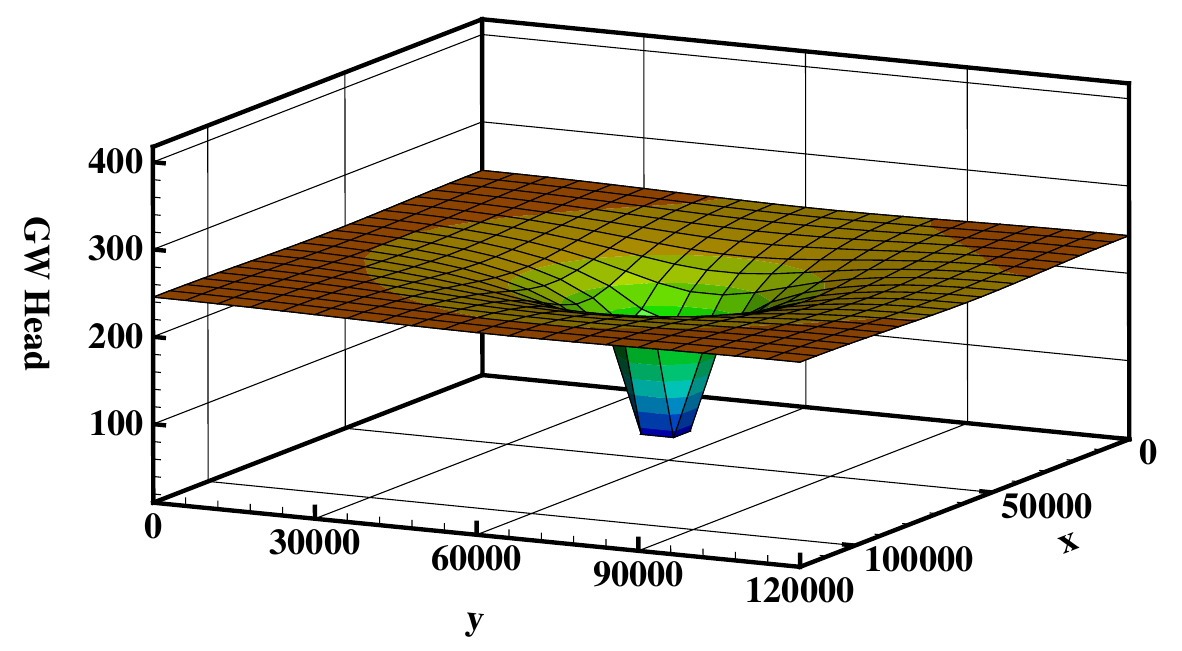}
    \end{subfigure}\hfil
   \begin{subfigure}{0.5\textwidth}
        \includegraphics[width=1\linewidth]{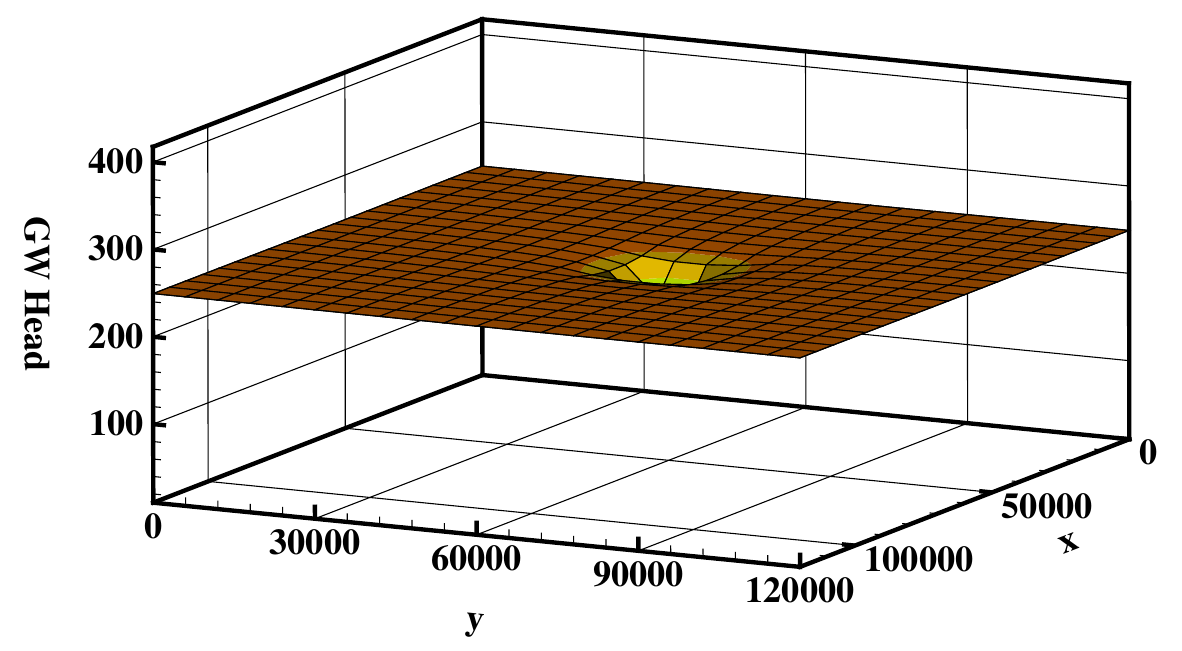}  
	\end{subfigure}\hfil
	\begin{subfigure}{0.5\textwidth}
	\includegraphics[width=1\linewidth]{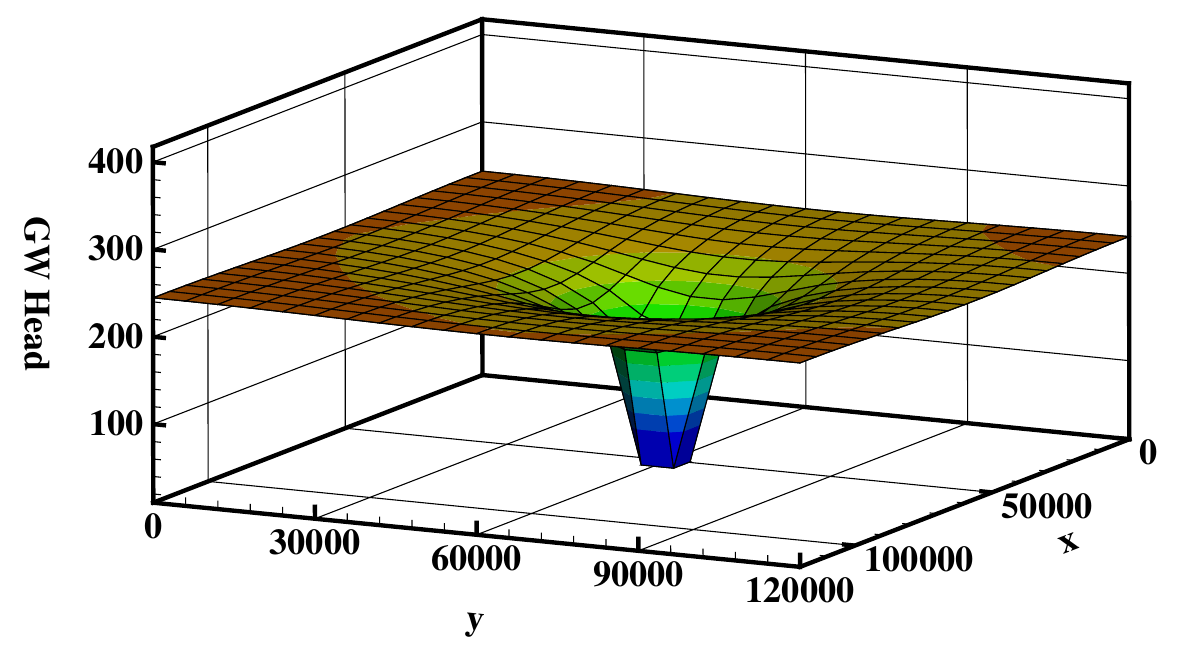}
    \end{subfigure}\hfil
   \begin{subfigure}{0.5\textwidth}
        \includegraphics[width=1\linewidth]{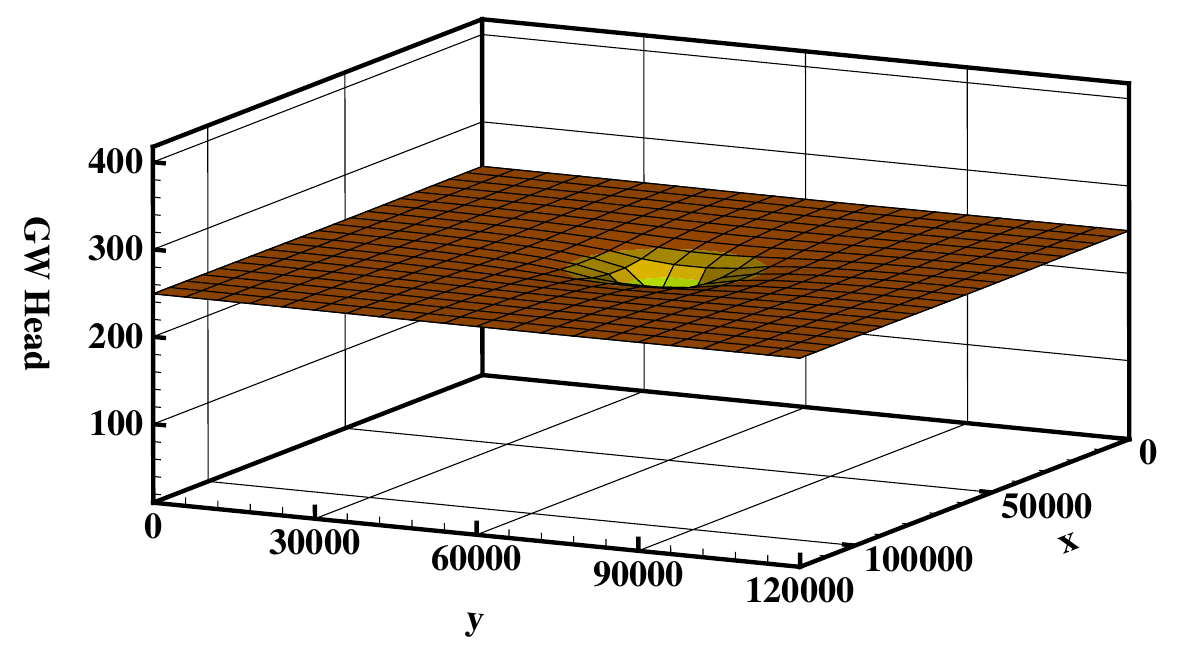}  
	\end{subfigure}\hfil
	\begin{subfigure}{0.5\textwidth}
	\includegraphics[width=1\linewidth]{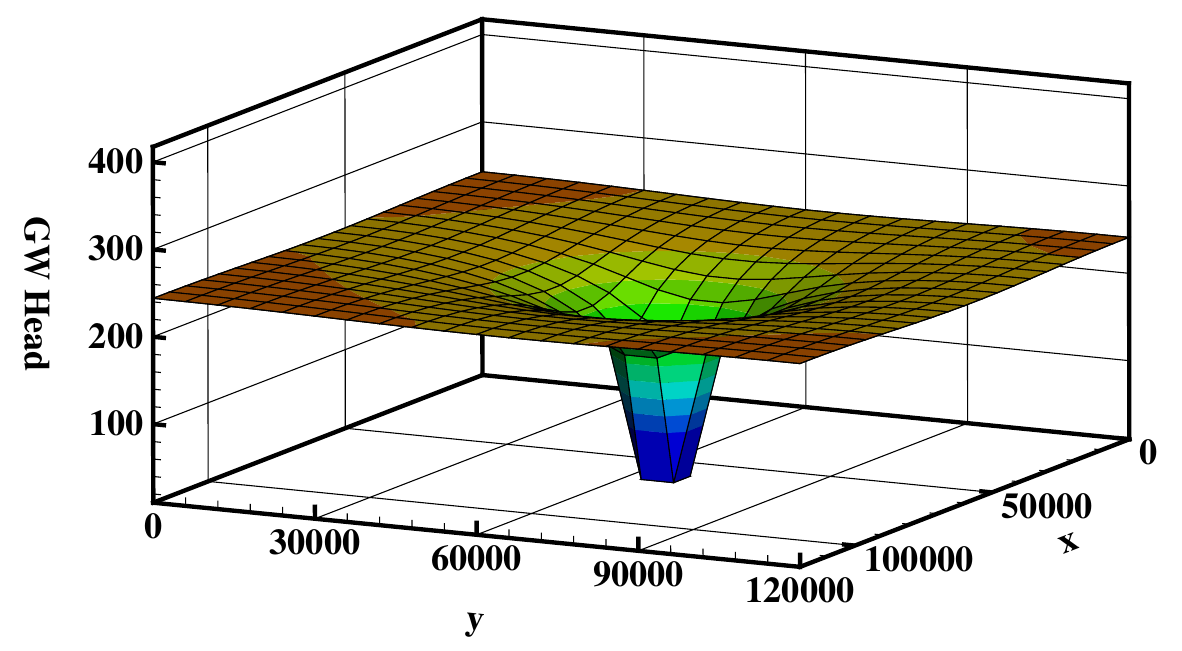}
    \end{subfigure}\hfil
\caption{{\sl Surface plot of groundwater head solutions for top (left) and bottom (right) layers using the JFNK approach in the two-layer system at the end of the $1^{st}$, $2^{nd}$, $3^{rd}$, and $4^{th}$ years (from top to bottom), respectively for the test case $2$. }}\label{fig:GW_2multi_1_4}
\end{figure} 


Figures \ref{fig:GW_2multi_1_4} display the comparison of the groundwater head levels calculated using the JFNK method for the top layer (plots on the left) and the bottom layer (plots on the right) at the end of each simulation year, exhibited from top to bottom. The largest draw-down is observed at the nodes where pumping occurs and it decreases gradually from the center towards the edges of the model domain. 

The top aquifer layer is an unconfined aquifer whereas the bottom layer is confined and has a lower storage coefficient compared to the top layer. This results in higher drawdown in the confined aquifer than the unconfined one even though the same amount of pumping is applied to both layers. Additionally, the model adheres to a physical constraint where the heads cannot fall below the bottom of each of the aquifer layers and pumping must be limited by the amount of storage available at the corresponding nodes (see equation \eqref{iwfm_l9}). This is apparent for the top layer in Figures \ref{fig:GW_2multi_1_4} . The bottom elevation of the top layer is at $200$ ft. It appears that groundwater heads in the top layer at nodes where pumping occurs have already dropped to $200$ ft at the end of the first year. From this point on, the model cuts down pumping and the heads at the top layer stay unchanged throughout the rest of the simulation period. On the other hand, lower aquifer layer has a larger storage available and the heads continue dropping throughout the simulation period.

\begin{figure}[!t]  
\centering
(a) \includegraphics[scale=.55]{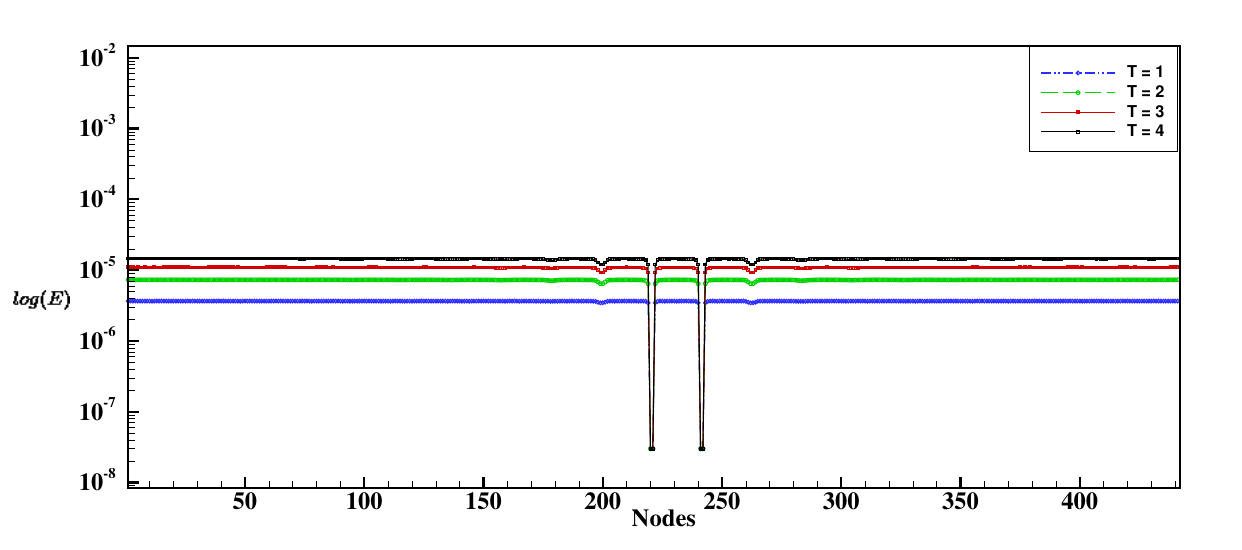} \\ 
(b) \includegraphics[scale=.55]{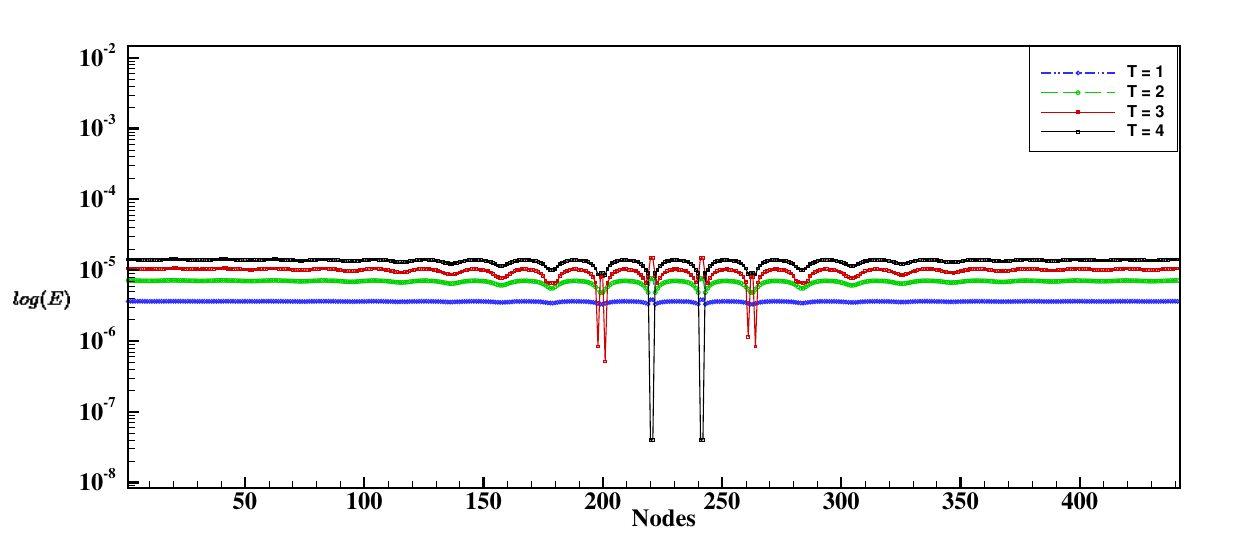}
 \caption{{\sl The log plot of the error between JFNK and NK solutions at $441$ nodes for top (a), and bottom (b) layers at the end of time $T = 1, 2, 3,$ and $4$ years for the test case $2$. }}\label{fig:Rel_Error_GW_2multi_1_4}
\end{figure}

\begin{figure}[!ht]
	\centering
	\includegraphics[width=0.6\linewidth]{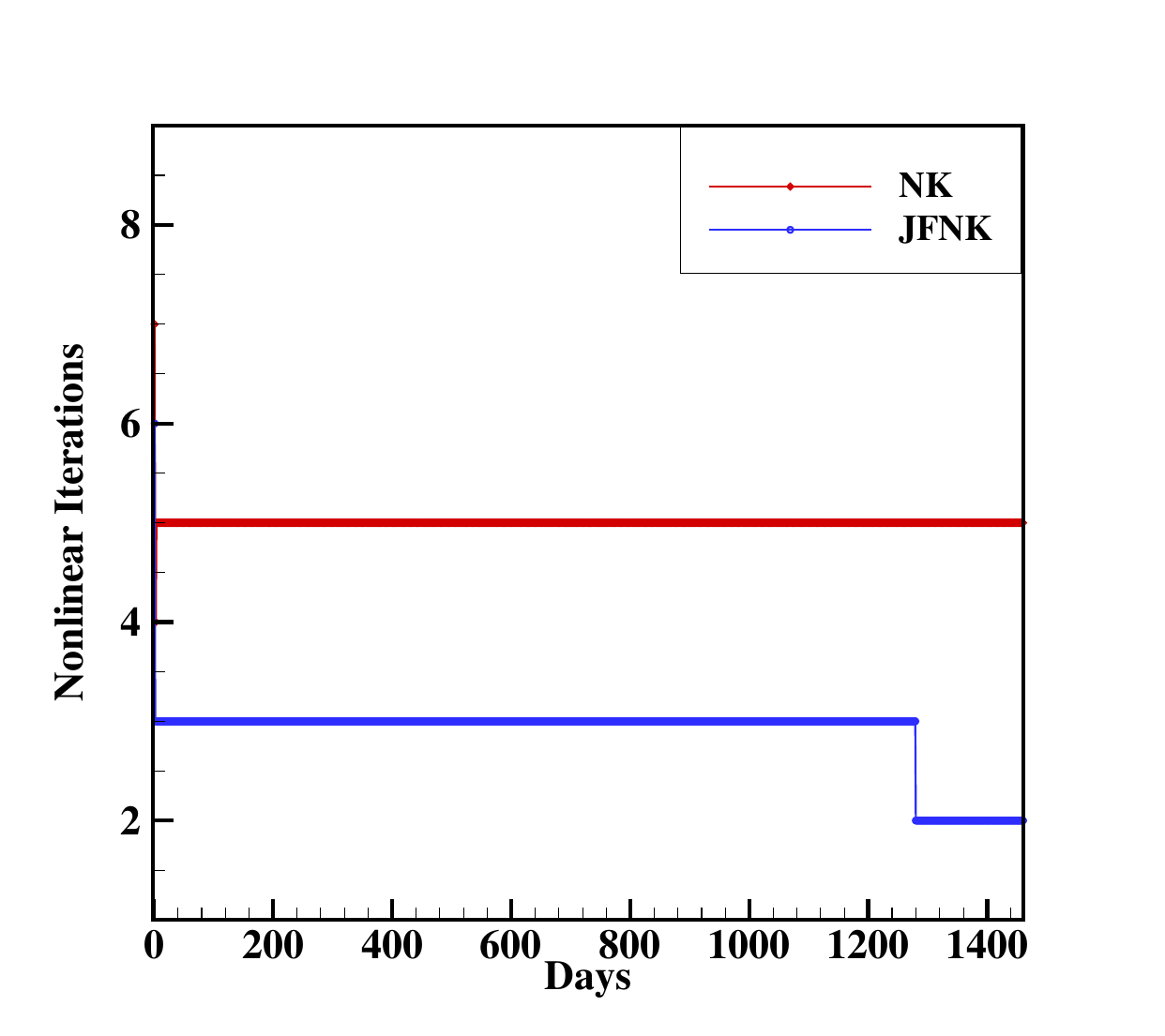}  
\caption{{\sl  The number of nonlinear iterations for each day in test case $2$ for both NK and JFNK approaches.} } \label{fig::2layer_Nonlinear_iterations::convg_hist}
\end{figure}

In Figure \ref{fig:Rel_Error_GW_2multi_1_4} (a) - (b), we show the error on a log scale between NK and JFNK for each node in the top and bottom layers at the end of $T = 1, 2, 3$, and $4$ years. The figure indicates that JFNK results are consistent with NK results at the end of each year.  The smallest error is occurring at the nodes where the pumping is happening because the head at the top layer hits the bottom of the aquifer where the pumping occurs. This threshold is a hard limit, and the head can no longer fall below this level. Thus, the error is the smallest at these nodes. One would assume that with a tighter convergence criteria, the results would be closer to the exact solution. The nonlinear iterations used by NK and JFNK methods for each day are displayed in Figure \ref{fig::2layer_Nonlinear_iterations::convg_hist}. The figure indicates that the JFNK method requires less outer iterations to achieve convergence at the given tolerance. 

\clearpage
\section{Performance of JFNK in IWFM}\label{sec:Implementatin_IWFM}\label{sec:Implementation_IWFM}
In this section, we will further examine the correctness and efficiency of the JFNK algorithm as implemented in IWFM. There are two models that we take into considerations. The first model is the same model that we mentioned in $\S$\ref{subsection:test_case1:GW}. The second model is the C2VSim model \cite{C2VSimFG_2020}, which is a realistic model with a greater number of degrees of freedom and is more advanced in nature.

\subsection{Groundwater flow model for an unconfined aquifer}
\begin{figure}[!t]  
\centering
\includegraphics[scale=0.55]{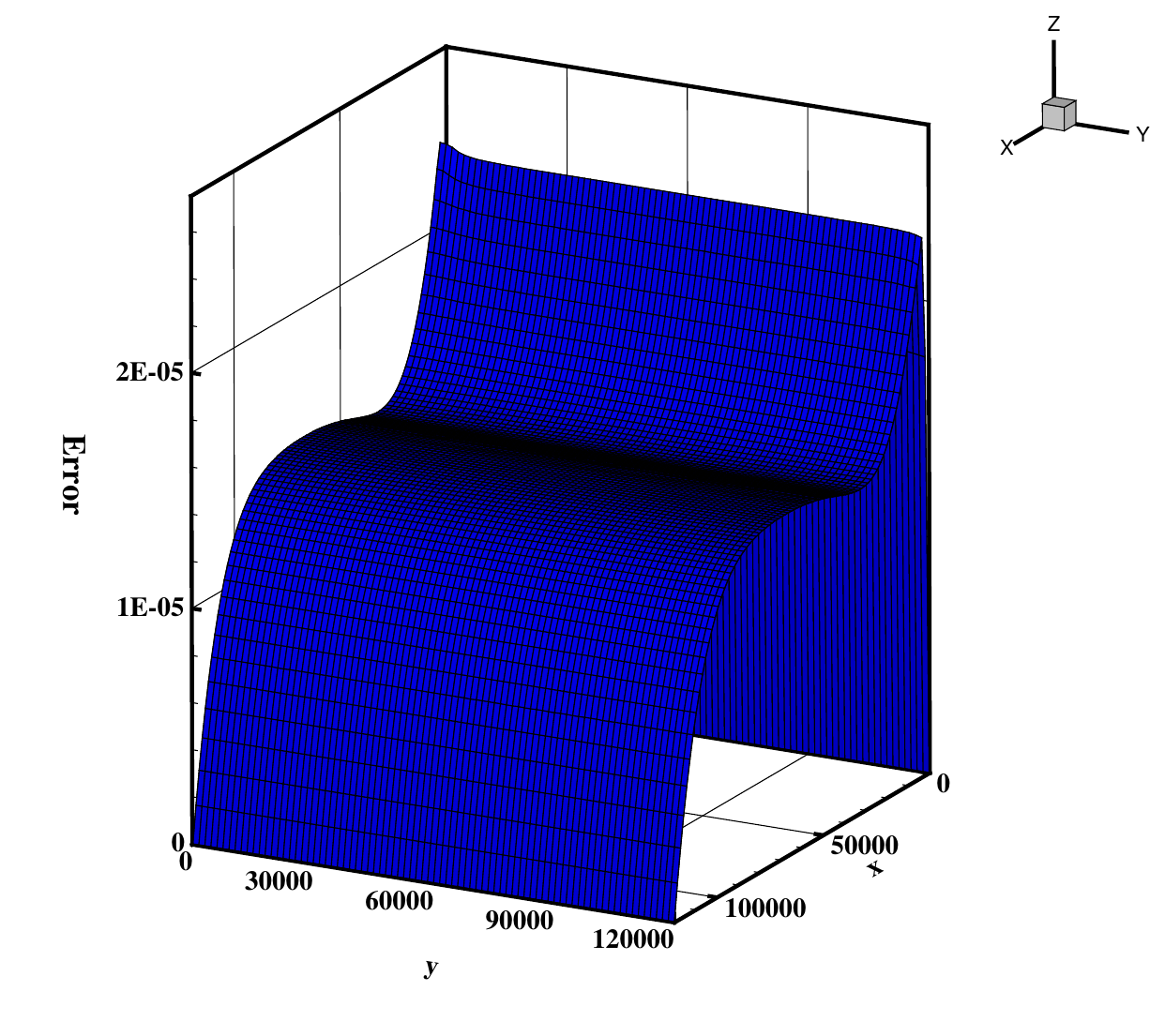}
 \caption{{\sl The error surface plot between JFNK and NK solutions at the end of $4$ years.}}\label{fig::iwfm_relative_error_plain_gw_case}
\end{figure}

\begin{figure}[!ht]
	\centering
	\begin{subfigure}{0.5\textwidth}
		\includegraphics[width=1.1\linewidth]{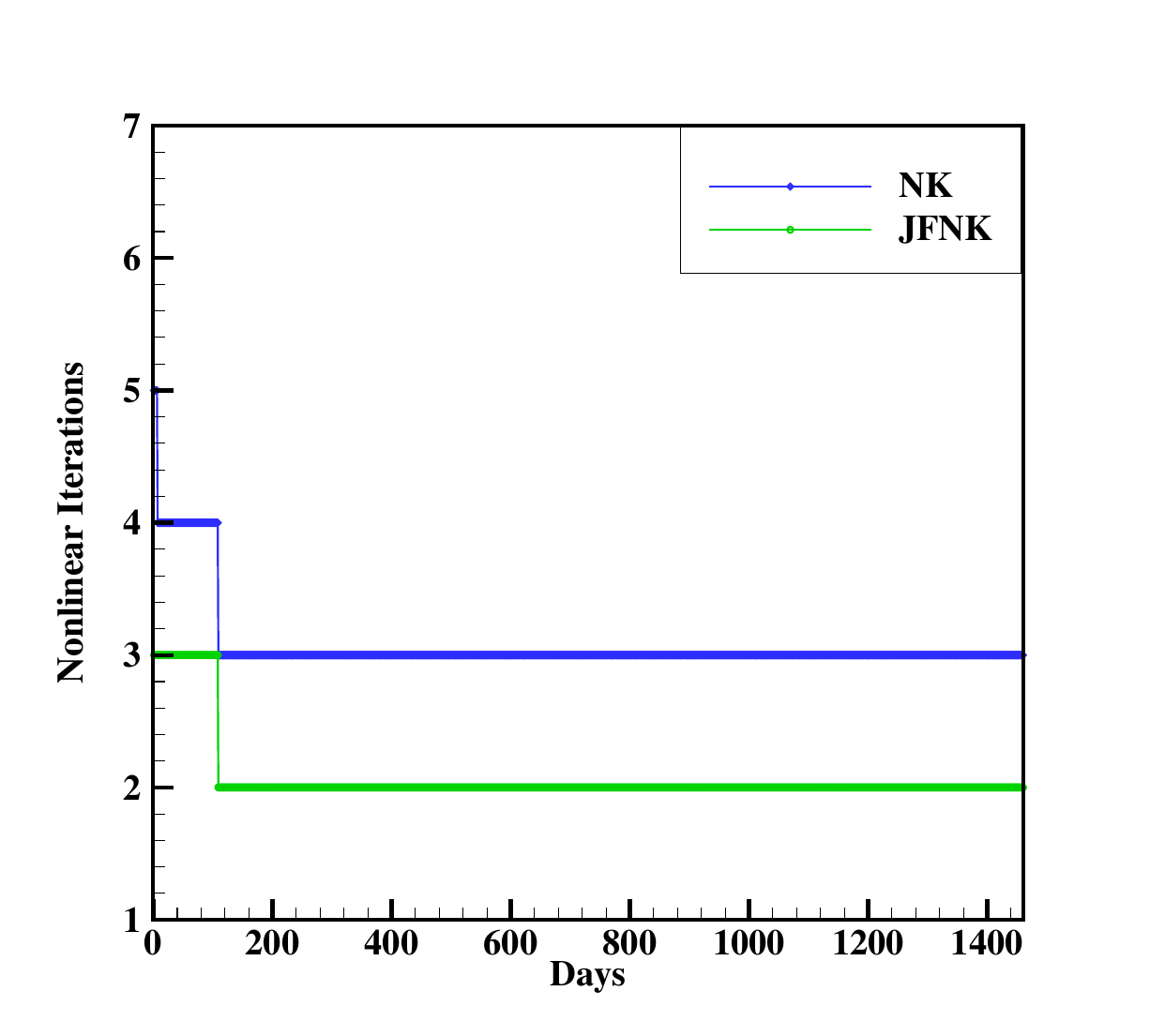}  
		\caption{}
	\end{subfigure}\hfil
	\begin{subfigure}{0.5\textwidth}
		\includegraphics[width=1.1\linewidth]{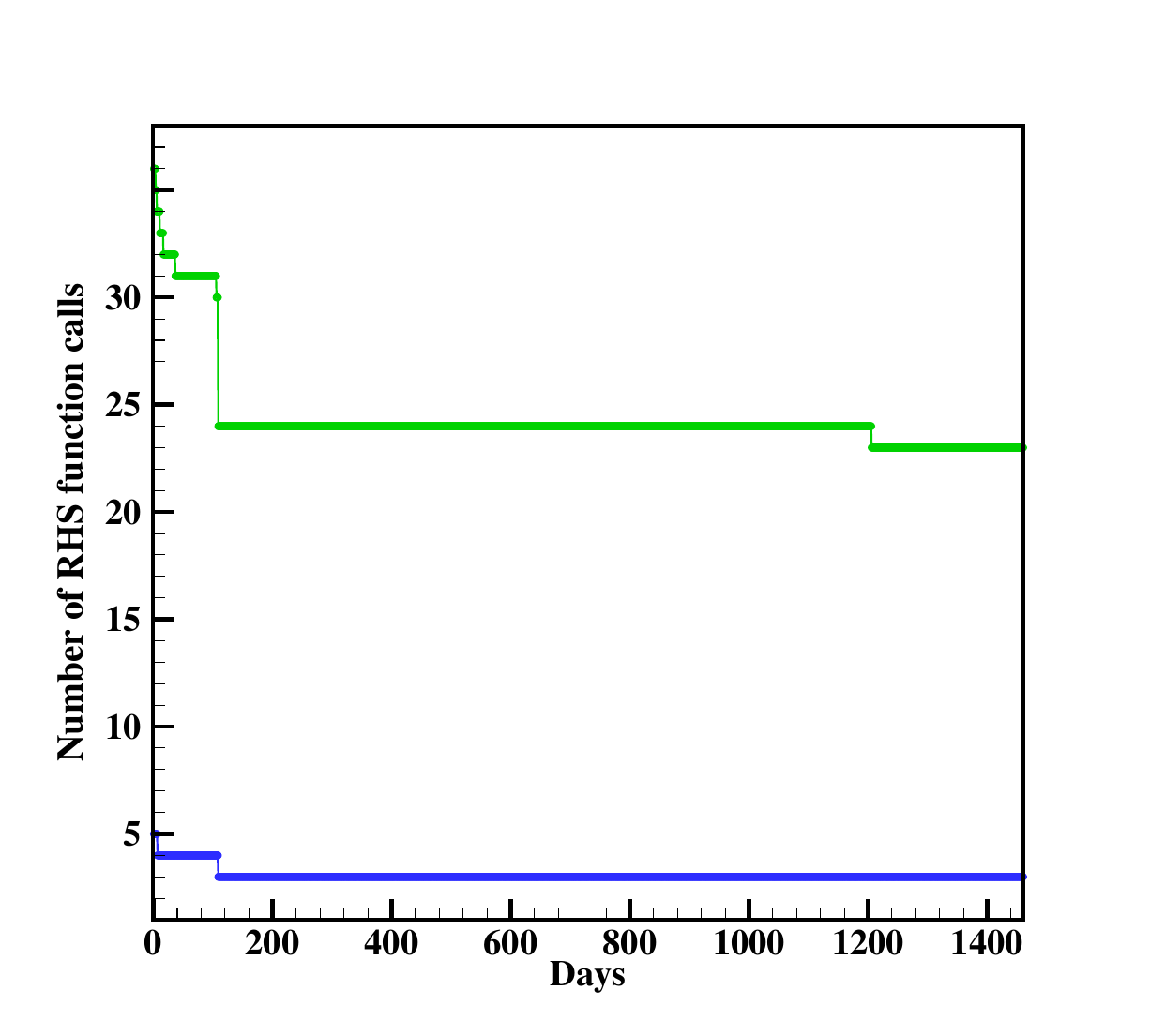}
		\caption{}
	\end{subfigure}\hfil
\caption{{\sl (a) The number of nonlinear iterations and (b) the frequency of calling the RHS function from both the NK and JFNK methods corresponding to each day. } } \label{fig::iwfm_NT_plain_gw_case}
\end{figure}

We considered the same model and parameters described earlier in $\S$\ref{subsection:test_case1:GW} to check the runtimes for the JFNK method implemented in IWFM. Figure \ref{fig::iwfm_relative_error_plain_gw_case} shows the surface plot of the absolute error between the NK and JFNK solutions at the end of four years, and the error magnitude ranges from $10^{-4}$ to $10^{-6}$. This proves the correctness of the JFNK method in IWFM. The nonlinear iterations and the total number of RHS function calls used by the NK and JFNK methods for each day are shown in figures \ref{fig::iwfm_NT_plain_gw_case} (a) and (b), respectively. The behavior is consistent with that observed in $\S$\ref{sec:Numerical examples}. Note that IWFM uses the ILU preconditioner \cite{saad1994ilut} for GMRES within the NK method. 

In IWFM, the overall computation times for NK and JFNK methods are $25.119$ and $47.93$ seconds, respectively. The total computing time in both techniques have increased compared to the times observed in $\S$\ref{subsection:test_case1:GW}. This is mainly due to the fact that RHS function in IWFM takes longer to execute than the one developed in $\S$\ref{subsection:test_case1:GW}. It includes many checks and balances to make sure that physics is not violated (e.g. groundwater heads do not fall below the bottom of the aquifer layer, pumping is scaled down properly with respect to available storage, stream flows do not become negative values, etc.). Since IWFM is a finite element model, forming the RHS vector is much more time consuming due to the unstructured nature of finite element grids and how the grid-related information is mapped onto the computer memory. Therefore, computation of the RHS function in IWFM takes longer than the one developed in $\S$\ref{subsection:test_case1:GW}.   The following are the key data
from the IWFM using NK and JFNK:

 \begin{enumerate}
     \item The total number of calls to the RHS function in the NK method is $4497$, while the total number in the JFNK method is $35631$. The factor of calling RHS function ($\varpi$) in JFNK method is \textbf{7.92} ($= 35631/4497$) times more compared to NK.  \label{tmp:c1}
     \item In JFNK, a single RHS function call consumes about $9.94 \times 10^{-4}$ seconds. \label{tmp:c3}
     \item In NK method, forming the coefficient matrix and the right hand side vector in \eqref{eq:nls4} consumes about $2.48 \times 10^{-3}$ seconds. \label{tmp:c4}
     \item  Based on items \ref{tmp:c3} and \ref{tmp:c4} above, compiling Jacobian one time in NK method is estimated to consume roughly $1.49 \times 10^{-3}$ seconds ($= 2.48 \times 10^ {-3} - 9.94 \times 10^{-4}$). \label{tmp:c5}
 \end{enumerate}
 Based on the above data, we can conclude that for the NK method, the total time taken by the RHS function is $4497 \times (9.94 \times 10 ^{-4}) =  4.470018$ seconds and the number of Jacobian calls in NK is $4497$. Therefore, the total time taken by Jacobian function is $4497 \times (1.49 \times 10^{-3}) =  6.7563$ seconds.
It suggests that the NK method spends approximately $50 \%$ of the overall runtime 
of $25.119$ seconds on the coefficient matrix and the right-hand side vector of \eqref{eq:nls4}. 
 
On the other hand, for the JFNK method, the total time taken by the RHS function is $35631 \times (9.94 \times 10 ^{-4}) =  35.41$ seconds, which accounts for approximately $74\%$ of the total runtime of $47.93$ seconds. In addition, we note that even though the RHS frequency factor $\varpi$ is similar to the one in \S\ref{subsection: eff}, but the total runtime increased by around 95\% compared to the NK method. Thus, it is important to minimize the runtime that RHS function takes to execute for the JFNK method to be faster than the NK method.

\subsection{C2VSim Model}
California's Central Valley (CV) plays a crucial role in the state's economy and in its water resources.  CV and its watersheds drain $38 \%$ of California's area and provide water to roughly $27$ million people in areas of the Central Coast and Southern California \cite{dwr2020}. The Sacramento-San Joaquin Delta and San Francisco Bay, the west coast's biggest estuary, obtain freshwater from CV runoff \cite{usbr2016}. Additionally, CV is  one of the world's leading agricultural regions and it is crucial to the U.S. food security as it produces most U.S. vegetables, fruits, and nuts \cite{usgs2020centralvalley}. 

California and CV hydrology are highly variable geographically and inter-annually. The state gets most of its precipitation during winter months in the form of snow in the Sierra Nevada mountains to the east of the CV and as rainfall within the CV. Majority of the precipitation occurs in the northern half of the CV while majority of the water demand is in summer months in the southern half. Federal, state and local officials have built large reservoirs to store the snow melt from the Sierra Nevada mountains to be used during the dry summer months to meet the water demand and large canals to carry this water from north to south. Despite this massive infrastructure to store and transmit water within the CV, surface water supply typically falls short of the demand, and groundwater becomes a substantial water supply. On average, groundwater  provides $40\%$  of the water demand in CV and this fraction can go up to $60\%$ in dry years. This has resulted in substantial groundwater development in the CV, and it is one of the most extensively pumped aquifer systems in the US.

\begin{figure}[!ht]
\begin{center}
\includegraphics[width=1\textwidth]{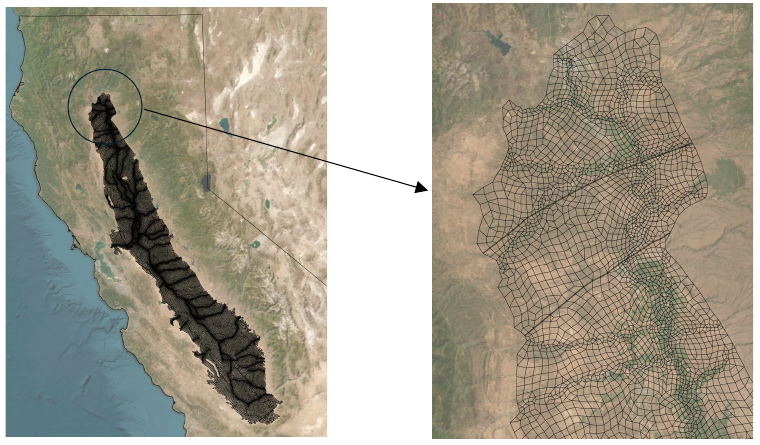} 
\end{center}
\caption{Present the geographic location of the study area in the C2VSim model and the mesh layout along with an enlarged view of the black-circle region.} \label{fig::c2vsim_mesh}
\end{figure}
To aid in their CV water resources management planning studies, California Department of Water Resources (CADWR) has developed the California Central Valley Groundwater-Surface Water Simulation Model (C2VSim) using IWFM as the underlying numerical engine \cite{C2VSimFG_2020}. C2VSim discretizes the CV model domain, an area of $20,750$ square miles, into $32,537$ finite element cells with maximum, minimum and median cell areas of 1770.8, 4.0 and 326.0 acres, respectively. The mixed quadrilateral-triangular finite element grid is shown in Figure \ref{fig::c2vsim_mesh}, along with an enlarged view of the circled region. The model represents the stream network with $110$ stream reaches discretized into $4,634$ stream nodes. The groundwater aquifer, whose thickness ranges from $200$ feet to $9000$ feet, is divided into four computational layers. The maximum, minimum and median layer thickness, considering all four layers, are $5576.1$, $41.9$ and $358.3$ feet, respectively. C2VSim is an integrated hydrologic model. Aside from groundwater and stream flows, it also simulates land surface and root zone flow processes, agricultural and urban water demands as well as groundwater pumping and surface water diversions to meet these demands. The values of $\tau_{h} =10^{-2}$, $p = 400$, $M_k = 1500$, and GMRES tolerance $=10^{-4}$ are considered in the simulation. C2VSim is calibrated for the simulation period starting from October $1973$ and ending in September $2015$ using a monthly timestep.

\begin{figure}[!ht]
	\centering
	\begin{subfigure}{0.5\textwidth}
		\includegraphics[width=1.1\linewidth]{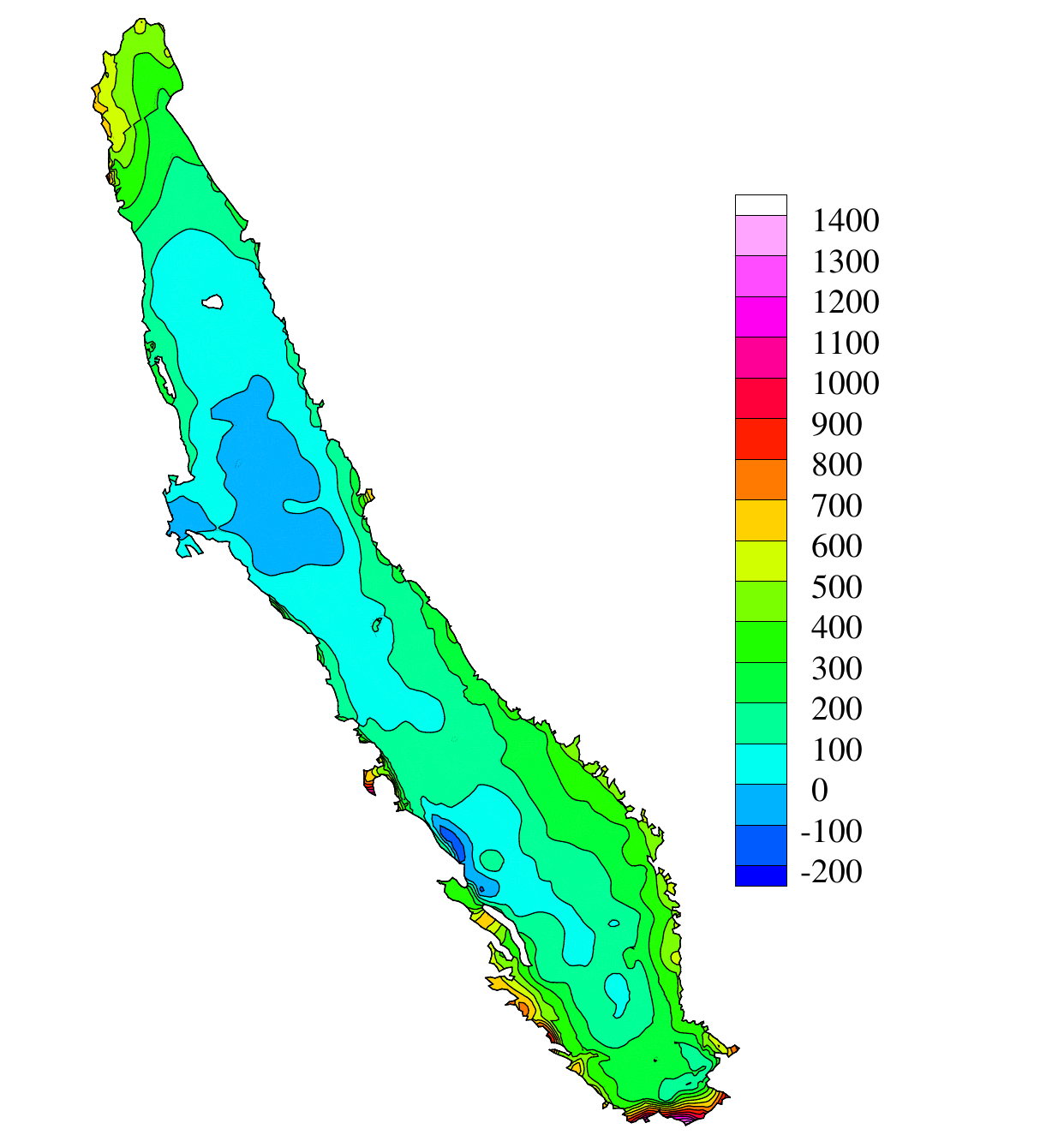}  
		\caption{}
	\end{subfigure}\hfil
	\begin{subfigure}{0.5\textwidth}
		\includegraphics[width=1.1\linewidth]{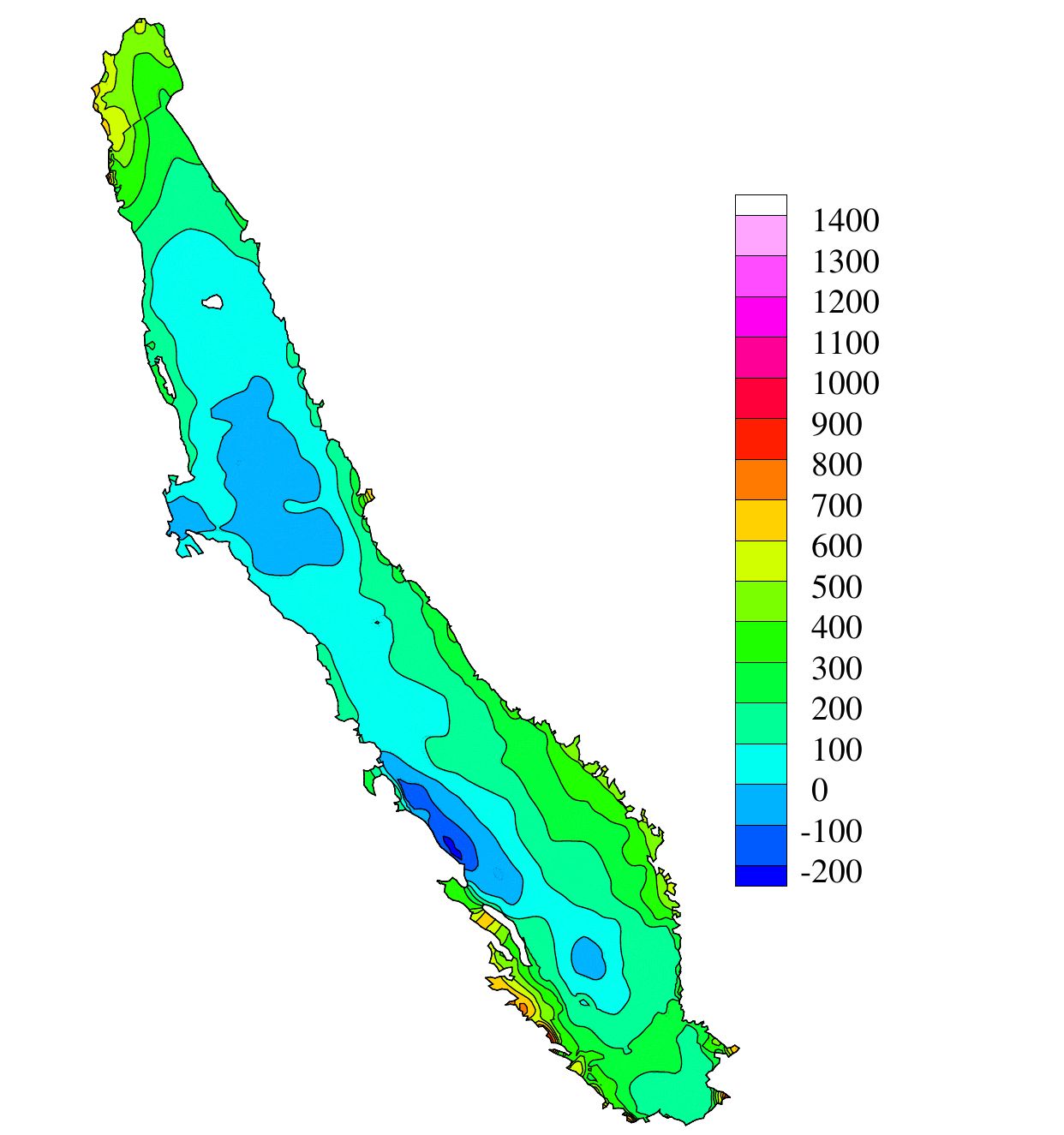}
		\caption{}
	\end{subfigure}\hfil
	\begin{subfigure}{0.5\textwidth}
		\includegraphics[width=1.1\linewidth]{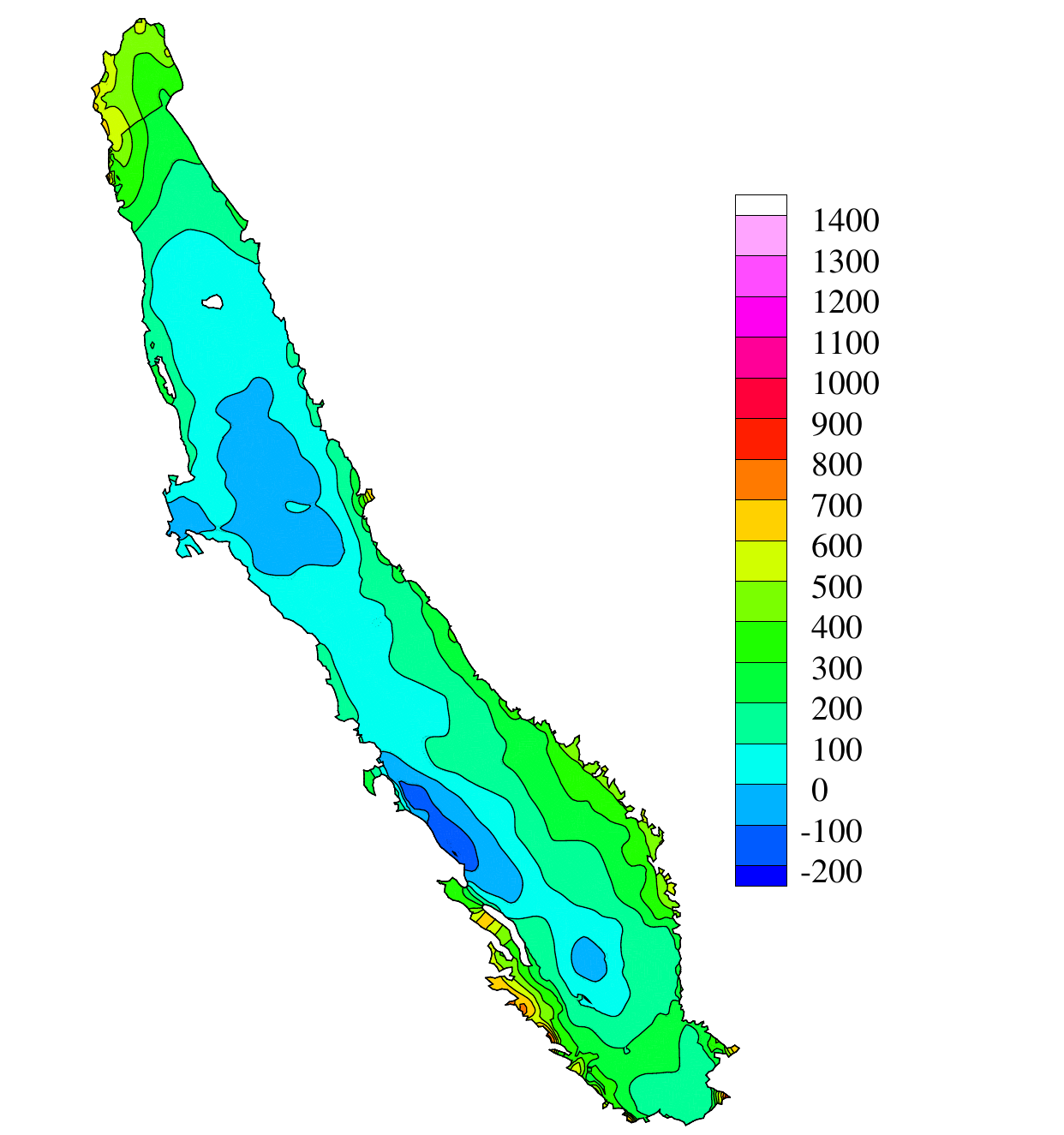}  
		\caption{}
	\end{subfigure}\hfil
	\begin{subfigure}{0.5\textwidth}
		\includegraphics[width=1.1\linewidth]{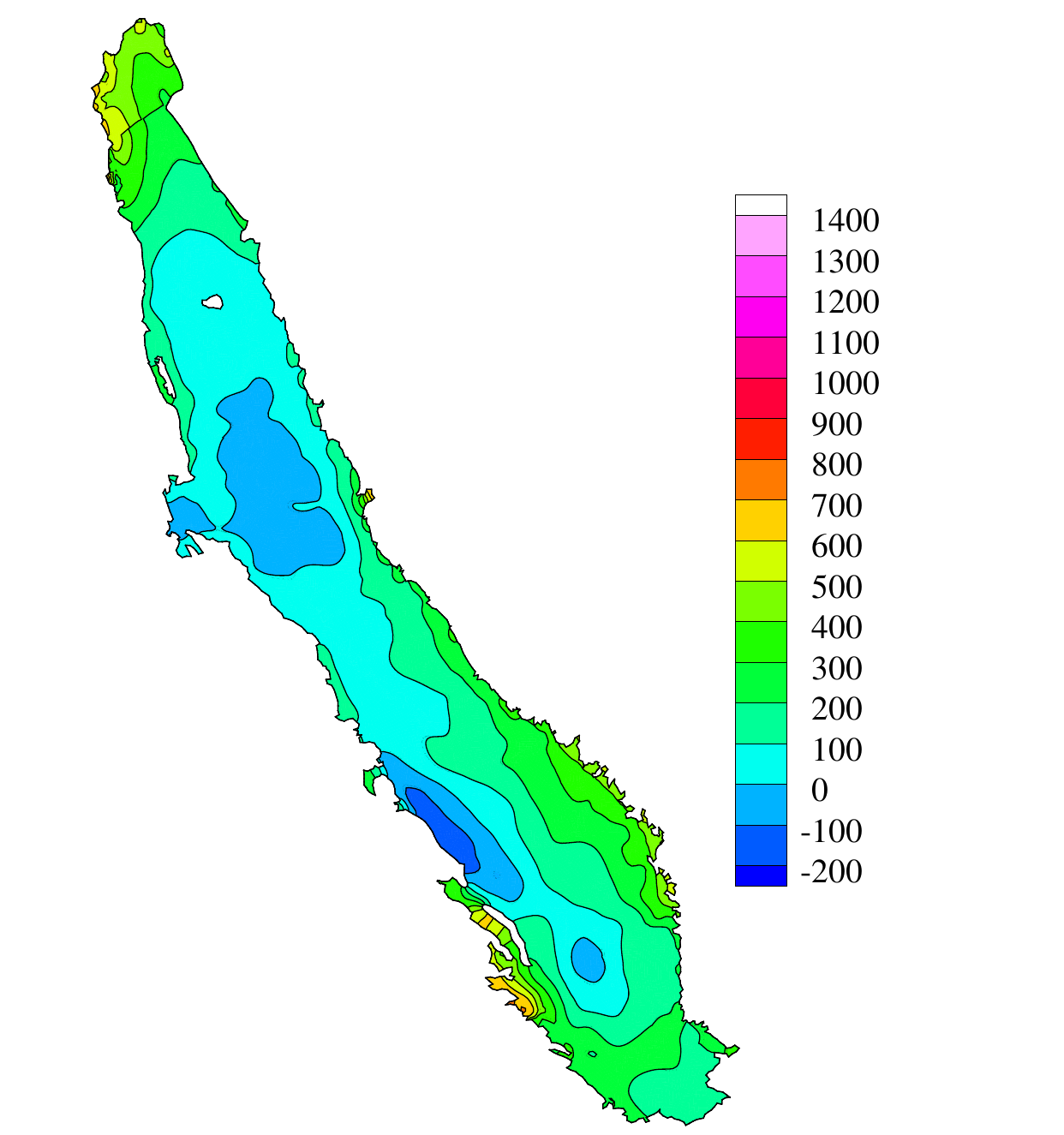}
		\caption{}
	\end{subfigure}\hfil
     \caption{{\sl Groundwater head solutions from JFNK method for layers $(a) 1^{st}$, $(b) 2^{nd}$, $(c) 3^{rd}$, and $(d) 4^{th}$ at the end of $1$ year of simulation for C2VSim model. }} \label{fig::GW_C2VSim}
\end{figure} 

As an initial test of the JFNK method, C2VSim was simplified by turning off all the simulation features except the groundwater component.  Figure \ref{fig::GW_C2VSim} shows the simulated groundwater heads at each aquifer layer after $1$ year. In Figure \ref{fig::GW_C2VSim} $1^{st}$ layer is the topmost aquifer layer, and the $4^{th}$ layer is the deepest aquifer layer. 

Figure \ref{fig::GW_C2VSim} shows that for all layers, groundwater head is lower in the center of the Valley and higher towards the edges. This is expected due to the geology of the CV: The aquifer is deeper towards of the center of the valley and it gets shallower towards the edges. Additionally, Figure \ref{fig::GW_C2VSim} shows that  groundwater heads in layers $2-4$ are very similar but they are different than those in layer $1$. The lowest groundwater heads in all layers occur in the southern half of the Valley. Although this is consistent with the observations, it should be re-emphasized that this version of C2VSim was not a representation of the real hydrologic system but a groundwater-only model for testing of JFNK model only.

\begin{figure}[!ht]
	\centering
	\begin{subfigure}{0.5\textwidth}
\includegraphics[width=1\linewidth]{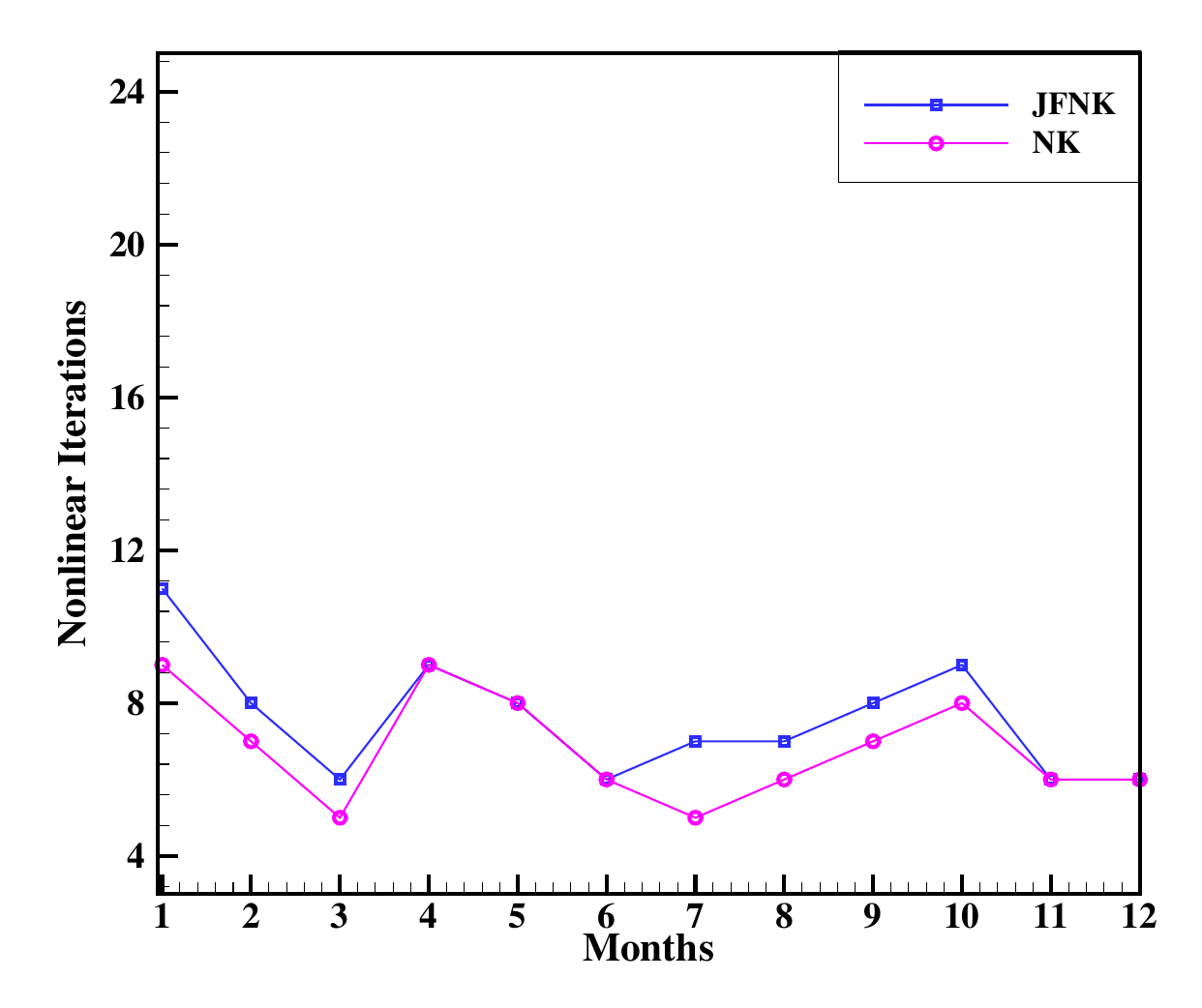} 
		\caption{}
	\end{subfigure}\hfil
	\begin{subfigure}{0.5\textwidth}
		\includegraphics[width=1\linewidth]{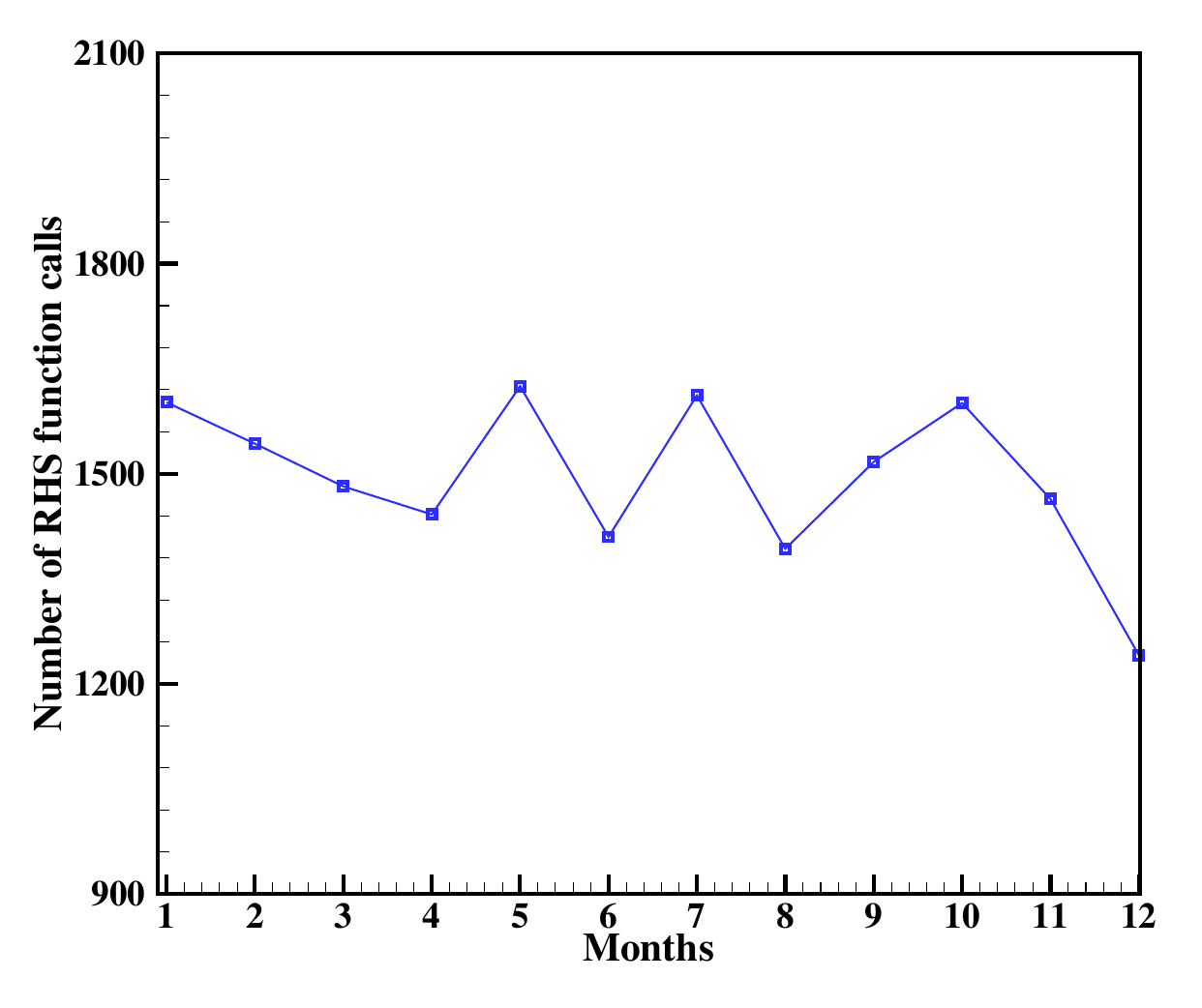}
		\caption{}
	\end{subfigure}\hfil
\caption{{\sl (a) The number of nonlinear iterations from NK and JFNK methods and (b) the number of RHS function call in the JFNK method corresponding to each month of the C2VSim model with a one-year simulation period.} } \label{fig::iwfm_NT_C2VSim_case}
\end{figure} 

Figure \ref{fig::iwfm_NT_C2VSim_case}(a) presents the number of nonlinear iterations corresponding to each month of the NK and JFNK methods in a one-year simulation period. The pattern of the nonlinear iteration from both methods is similar, and it is also important to note that IWFM uses the ILU preconditioner within the NK method, whereas the JFNK method does not. Figure \ref{fig::iwfm_NT_C2VSim_case}(b) shows the frequency of the RHS function in the JFNK method corresponding to each month. Note that the number of RHS function calls in the NK method is equivalent to the number of nonlinear iterations as shown in Figure \ref{fig::iwfm_NT_C2VSim_case}(a). The increased number of RHS function calls for the JFNK method is expected and this behavior is also seen in the previous examples.

\begin{figure}[!ht]
\begin{center}
\includegraphics[width=1\textwidth]{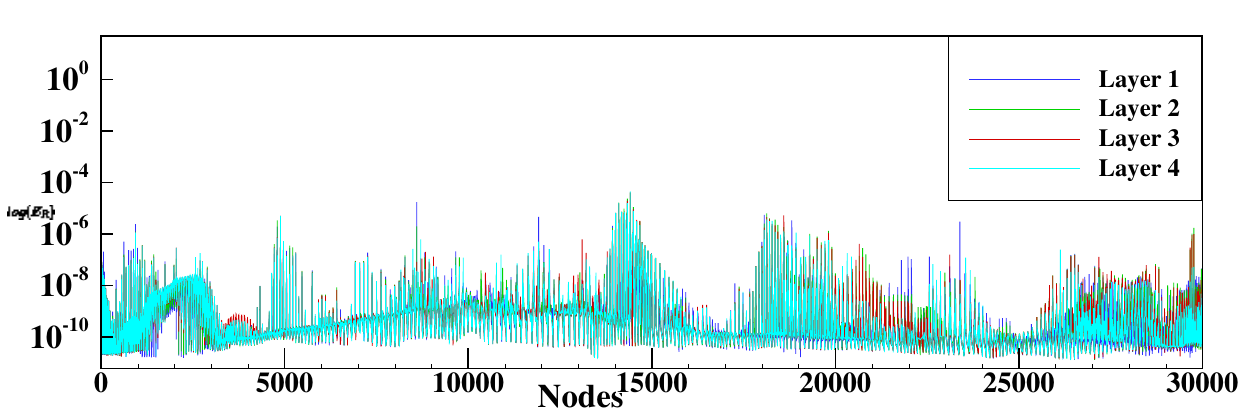} 
\caption{{\sl The log plot of the relative error on each node between JFNK and NK solutions for all aquifer layers of the C2VSim model at the end of the simulation period (layer 1 is the top and layer 4 is the deepest aquifer layers) .}}\label{fig::rel_error_c2vsim}
\end{center} 
\end{figure}

To compare the groundwater heads simulated with the NK and JFNK methods, we define a relative error as 
\begin{equation}
E_{R} = \frac{|h_{NK}(x,y,t) - h_{JFNK}(x,y,t)|}{|h_{NK}(x,y,t)|} \label{eqn::rel_error}
\end{equation}
In Figure \ref{fig::rel_error_c2vsim}, we present the log plot of relative error comparing the results of NK and JFNK for all $32,537$ nodes of all four layers. The magnitude of the error lies in the range of $O(10^{-4})$ to $O( 10^{-11})$, which indicates that there is a good match between the results obtained by NK and JFNK. It is important to note the fact that the JFNK implementation in IWFM ran on a complex model setup and provided results for the analysis itself, which can be considered a success.

\clearpage
\section{Conclusion}\label{sec:Conclusion}
Complex physical problems can be modelled mathematically through nonlinear systems of equations. The standard Newton method has been the most frequently used approach to solving the resulting system over a period of a century. Adopting a suitable iterative method to solve a nonlinear system is of utmost importance for obtaining the desired accuracy and optimal efficiency. This study examines the iterative technique known as the "Jacobian-Free Newton-Krylov method" for addressing the groundwater flow problem with multiple aquifers.

Generally, the selection of an algorithm that optimizes performance in terms of accuracy and efficiency is highly dependent upon the question of interest. It is not always feasible to compute the exact or true Jacobian in Newton method at each nonlinear iteration, particularly when the system is complex and high-dimensional, as one may overlook system features when calculating the partial derivatives. In contrast to Newton method, the primary characteristic of the JFNK algorithm is the elimination of the necessity to generate, preserve, and invert the extensive Jacobian matrix during each nonlinear iteration, making this method distinctly unique. The key component in JFNK is that the Jacobian-vector product is being approximated through the forward or central difference, which utilizes the evaluation of $F$.  
 
We conducted a comprehensive examination and performed numerical simulations on several test cases (one of which included discontinuities) using IWFM, a finite element integrated hydrologic model. Our focus is on the groundwater flow equation in single and multiple layers within an unconfined aquifer. The results exhibit the correctness and robustness of the algorithm for these models and are reasonable when compared to the results from the application of the standard Newton method using the exact Jacobian. 

We evaluate the efficiency of the iterative and direct approaches within a single-layer system lacking pumping parameters. To conduct a fair comparison of profiling data, we established the same problem programmed in Matlab and IWFM from both approaches over a span of $4$ years and discover the circumstances under which the performance of the JFNK algorithm can be enhanced and capitalized. Despite the RHS frequency factor ($\varpi$) indicating that the RHS subroutine is executed approximately 6.61 times more than the NK approach, the JFNK method still outperforms the optimal NK method by nearly $17\%$, particularly when the RHS subroutine is contained within a single file or is uncomplicated in nature during frequent calls. However, when we coded the identical problem in IWFM and implemented the JFNK algorithm in order to solve the nonlinear system, the behaviour of $\varpi$ remained consistent. Nevertheless, we experienced a $48 \%$ reduction in efficiency over an optimal NK method. This clearly shows how important it is that the RHS subroutine has been coded efficiently and in a parallelized way. This can be crucial for improving speed in real-life applications where the subroutine may be spread across multiple files or interconnected to other subroutines that have to deal with many physical factors.


Finally, we investigated the robustness of the JFNK approach with a real-life application for California's Central Valley with degrees of freedom on the order $10^6$. 

Although JFNK technique alleviates the need for explicit calculation and storage of the Jacobian matrix and improves CPU efficiency, in this study, we have not been able to prove that JFNK saves CPU runtime compared to NK method in IWFM. It may be worthwhile to explore the use of preconditioners in the future to further enhance the algorithm's performance. Considering Anderson acceleration \cite{walker2011anderson} as preconditioned appears advantageous due to the absence of matrix storage, and it would be interesting to explore how the $\varpi$ factor will behave.

\paragraph{Acknowledgment} This research was funded by the California Department of Water Resources and the University of California, Davis contract $4600014853$. 

\bibliographystyle{plain}

\end{document}